\documentclass[12pt,a4paper,reqno]{amsart}
\usepackage[english]{babel}
\usepackage{amssymb,latexsym,amsfonts,amsthm,upref,amsmath}
\usepackage[margin=1in]{geometry}
\usepackage{hyperref}
\hypersetup{colorlinks,citecolor=blue,filecolor=black,linkcolor=blue,urlcolor=blue}
\usepackage{comment} 
\usepackage{enumerate} 
\usepackage{tikz}
\usepackage{float}
\usepackage{cleveref}

\newtheorem*{corintro*}{Corollary}
\newtheorem*{thm*}{Theorem}
\newtheorem*{lem*}{Lemma}
\newtheoremstyle{prim}{}{}{\normalfont}{}{\bfseries}{.}{ }{}
\theoremstyle{prim}
\newtheorem{ex}{Example}
\newtheoremstyle{stil}{}{}{\slshape}{}{\bfseries}{.}{ }{}
\theoremstyle{stil}
\newtheorem{thm}{Theorem}[section]
\newtheoremstyle{defi}{}{}{}{}{\bfseries}{.}{ }{}
\theoremstyle{defi}
\newtheorem{defn}[thm]{Definition}
\theoremstyle{defi}
\newtheorem{rem}[thm]{Remark}
\theoremstyle{stil}
\newtheorem{pro}[thm]{Proposition}
\theoremstyle{stil}
\newtheorem{lem}[thm]{Lemma}
\theoremstyle{stil}
\newtheorem{kor}[thm]{Corollary}
\newenvironment{prf}{\noindent \textit{Proof.}}{\null\hfill$\qed$\hskip
2mm\vskip 2mm}

\newcommand{\wt}{\mathop{\mathrm{wt}}}

\newcommand{\ch}{\mathop{\mathrm{ch}}}

\newcommand{\xpan}{\mathop{\mathrm{span}}}
\newcommand{\q}{{\protect\underline{\mathsf{q}}}}
\newcommand{\vac}{\mathop{\mathrm{\boldsymbol{1}}}}
\newcommand{\ndo}{\mathop{\mathrm{End}}}
\newcommand{\om}{\mathop{\mathrm{Hom}}}

\newcommand{\im}{\mathop{\mathrm{Im}}}

\newcommand{\sym}{\mathop{\mathrm{Sym}}}
\newcommand{\gauss}[2]{\genfrac{[}{]}{0pt}{}{#1}{#2}}
\newcommand{\diag}{\mathop{\mathrm{diag}}}

\numberwithin{equation}{section}

\begin{document}

\title[Higher level vertex operators for $U_q (\widehat{\mathfrak{sl}}_2)$]{Higher level vertex operators for $U_q (\widehat{\mathfrak{sl}}_2)$}

\author{Slaven Ko\v{z}i\'{c}} 

\address{Department of Mathematics, University of Zagreb, 10000 Zagreb, Croatia}

\address{School of Mathematics and Statistics F07, University of Sydney, NSW 2006, Australia}

\email{kslaven@maths.usyd.edu.au}



\begin{abstract}
We study   graded nonlocal $\underline{\mathsf{q}}$-vertex algebras and we prove that they can be generated by certain sets of vertex operators. As an application, we consider the family of graded 
nonlocal $\underline{\mathsf{q}}$-vertex algebras $V_{c,1}$, $c\geq 1$, associated with the principal subspaces $W(c\Lambda_0)$ 
of the integrable highest weight $U_q (\widehat{\mathfrak{sl}}_2)$-modules $L(c\Lambda_0)$.
Using quantum integrability, we derive combinatorial bases for $V_{c,1}$ and compute the corresponding character formulae.
\end{abstract}

 \maketitle

   \section*{Introduction}
\allowdisplaybreaks
In their work \cite{LP}, J. Lepowsky and M. Primc found the, so-called integrability condition
for the affine Kac-Moody Lie algebra $\widehat{\mathfrak{sl}}_2$,
\begin{equation}\label{integrabilityrelation}
x_{\alpha} (z)^{c+1}=0
\end{equation}
on a level $c$ integrable    $\widehat{\mathfrak{sl}}_2$-module. 
In general, \eqref{integrabilityrelation} holds  on an arbitrary 
level $c$ integrable   $\widehat{\mathfrak{g}}$-module, when the simple root $\alpha$
is replaced by the maximal root of the untwisted affine Kac-Moody Lie algebra $\widehat{\mathfrak{g}}$.
Relation \eqref{integrabilityrelation} led to a construction of combinatorial bases for integrable highest weight 
$\widehat{\mathfrak{sl}}_2$-modules (cf. \cite{FKLMM},\cite{MP}) and, consequently, to a new series of combinatorial 
Rogers-Ramanujan-type identities.
Furthermore,  integrability relations played an important role in the construction of monomial bases
for certain substructures of $\widehat{\mathfrak{g}}$-modules such as principal subspaces (cf. \cite{FS},\cite{G},\cite{Bu1},\cite{Kaw})
and Feigin-Stoyanovsky's type subspaces  (cf. \cite{P1},\cite{P2},\cite{JP},\cite{T}).
For more information on principal subspaces the reader may consult, for example, the papers \cite{x1}--\cite{x4},\cite{S1},\cite{S2} and the references therein.
Using  Drinfeld realization  of quantum affine algebra $U_q (\widehat{\mathfrak{sl}}_{n+1})$ (see \cite{D})
and Frenkel-Jing realization of its integrable highest weight modules (see \cite{FJ})
J. Ding and T. Miwa found in \cite{DM} quantum integrability relations,
 \begin{equation}\label{dm7}
 x_{i}^{\pm}(z_{1})x_{i}^{\pm}(z_{2})\cdots x_{i}^{\pm}(z_{c+1})=0\quad\text{if}\quad z_{1}/z_{2}=z_{2}/z_{3}=\ldots=z_{c}/z_{c+1}=q^{\mp 2},
\end{equation}
on a level $c$ integrable   $U_q (\widehat{\mathfrak{sl}}_{n+1})$-module.

In this paper, we continue our research on vertex algebraic structures arising from Frenkel-Jing operators $x_{1}^{\pm}(z)$ for 
$U_q (\widehat{\mathfrak{sl}}_{2})$, which was initiated in \cite{Ko4}. 
So far  there were several fruitful approaches to associating   vertex algebra-like theories with the various quantum objects, such as quantum affine algebras or Yangians, which resulted in  some fundamental results and important constructions (cf. \cite{AB},\cite{Bor},\cite{EK},\cite{FR},\cite{Li1}--\cite{Li3}).   
However, 
motivated by the role of integrability \eqref{integrabilityrelation} in the representation theory of the affine Kac-Moody Lie algebras,
we introduce  graded nonlocal $\q$-vertex algebras, certain new structures 
 which are designed to make use of quantum integrability 
\eqref{dm7}.

The paper is organized as follows. In Section \ref{preliminaries}, we establish the notation and recall some well-known results from the theory of quantum affine algebras. Although the exposition is to a great extent similar
to the preliminary section in  \cite{Ko4}, we decided to include it in order to make this
paper  as self-contained as possible.

In Section \ref{gkva}, we introduce the notion of graded nonlocal $\q$-vertex algebra, which was motivated by \cite{Bak} and \cite{Li1}.
 Roughly speaking, graded nonlocal $\q$-vertex algebra is a triple
$(V,Y,\vac)$, where $V=\coprod_{u\in\mathbb{Z}}V_{(u)}$ is a vector space over field $\mathbb{F}\supset\mathbb{C}(\q)$ of characteristic zero, which satisfies all the axioms of vertex algebra (cf. \cite{LiLep}) except  Jacobi identity, which is replaced by
associativity,
\begin{equation}\label{assoc:78}
Y(a,z_0 +z_2)Y(b,z_2 \q^{u}) c=Y(Y(a,z_0)b,z_2)c\end{equation}
for all $a\in V_{(u)}$, $b,c\in V$, $u\in\mathbb{Z}$.
The variables $z_0$ and $z_2$ in \eqref{assoc:78} satisfy the noncommutative constraints
\begin{equation*}
z_2 z_0=\q z_0 z_2
\end{equation*}
 and certain natural grading restrictions are imposed on $V$ and $Y$. Naturally, we expect that in a more general setting a weaker form of \eqref{assoc:78} should be regarded.
\makeatletter
\def\thmhead@plain#1#2#3{%
  \thmname{#1}\thmnumber{\@ifnotempty{#1}{ }\@upn{#2}}%
  \thmnote{ {\the\thm@notefont#3}}}
\let\thmhead\thmhead@plain
\makeatother  

For a vector space $L$ over field $\mathbb{F}\supset\mathbb{C}(\q)$ of characteristic zero define $$\mathcal{E}(L)_t=\om(L,L((z)))\otimes\mathbb{F}[t].$$
Our goal is to construct graded nonlocal $\q$-vertex algebras generated by certain subsets of 
$\mathcal{E}(L)_t$.
The space $\mathcal{E}(L)_t$ is considered instead of $\mathcal{E}(L)=\om(L,L((z)))$ because it  allows us to efficiently introduce the $r$th products, $r\in\mathbb{Z}$, among vertex operators, which utilize  quantum integrability \eqref{dm7}. The products are defined for the operators satisfying a certain technical
requirement, quasi-commutativity, which may be considered as a more restrictive version of quasi-compatibility  introduced by H.-S. Li
(cf. \cite{Li1}). The main result in this section is the following theorem:
\begin{thm*}[\textbf{\ref{main2}}]
Let $\mathcal{S}$ be a quasi-commutative subset of $\mathcal{E}(L)_t$. There exists a unique
smallest graded nonlocal $\q$-vertex algebra $V\subseteq\mathcal{E}(L)_t$ such that $\mathcal{S}\subseteq V$.
Furthermore, if $\mathcal{S}$ consists of homogeneous elements, we have
\begin{equation*}
V=\xpan \left\{a_1(z,t)_{r_1}\cdots a_k(z,t)_{r_k}\vac\,:\,a_j(z,t)\in\mathcal{S},\,r_j< 0,\,j=1,...,k,\,k\in\mathbb{Z}_{\geq 0} \right\}.
\end{equation*}
\end{thm*}
As an application of Theorem \ref{main2} we construct graded nonlocal $\q$-vertex algebras 
generated by  quantum current operators $x_{i}^{+}(z)$ or $x_{i}^{-}(z)$ acting on an arbitrary restricted $U_q(\widehat{\mathfrak{g}})$-module, where $\widehat{\mathfrak{g}}$ is an affine Kac-Moody Lie algebra
 of type $(ADE)^{(1)}$. 
Some of the ideas and results in this section rely upon
Li's theory of nonlocal vertex algebras (cf. \cite{Li1},\cite{Li0}) even though different products and  structures are considered.

In Section \ref{section3}, we study graded nonlocal $\q$-vertex algebras $V_{c,1}$, $c\in\mathbb{Z}_{>0}$,
generated by the  operator $x(z,t):=x_{1}^{+}(z)\otimes t\in\mathcal{E}(L)_t$ for 
$U_q(\widehat{\mathfrak{sl}}_2)$, which acts on the integrable  
$U_q(\widehat{\mathfrak{sl}}_2)$-module $L=L(\Lambda_0)^{\otimes c}\supset L(c\Lambda_0) $ of level $c$.
Using (quantum) quasi-particles $x_m (z,t)$, $m=1,...,c$, from \cite{Ko1} we define the following subset
of $V_{c,1}$:
\begin{align*}
\mathcal{B}_{c,1}=\big\{&x_{m_1}(z,t)_{r_1}x_{m_2}(z,t)_{r_2}\ldots x_{m_k}(z,t)_{r_k}\vac\,:\big.\\
&\big. r_1,...,r_{k-1}\leq -2,\,r_k\leq -1,\,1\leq m_j\leq c,\,j=1,...,k,\,k\in\mathbb{Z}_{\geq 0} \big\},
\end{align*}
where $\vac$ is the vacuum vector in $V_{c,1}$.
By employing the second part of Theorem \ref{main2} and Koyama's realization of intertwining operators for $U_q(\widehat{\mathfrak{sl}}_2)$ (see \cite{Koyama}), we prove the main result of this section:
\begin{thm*}[\textbf{\ref{tralalala}}]
The set $\mathcal{B}_{c,1}$ forms a basis for $V_{c,1}$.
\end{thm*}
For $c=1$ the similar basis was already found in \cite{Ko4}. However, this was a basis 
for  $W_{1,\q}$, a certain subspace of a much bigger (nongraded) nonlocal $\q$-vertex algebra generated by $x(z)=x_{1}^{+}(z)\in\mathcal{E}(L(\Lambda_0))$ (which satisfied a slightly modified version of \eqref{assoc:78}).
The space $W_{1,\q}$ did not have any additional vertex algebraic structure, which was caused by a lack of grading restrictions. In this paper, as a consequence of the imposed grading restrictions, the operator $x(z)\otimes t$ (on level $1$) generates  graded nonlocal $\q$-vertex algebra
$V_{1,1}$, which gives rise to the same character formula  as $W_{1,\q}$.

It is important to emphasize that the use of the space $\mathcal{E}(L)_t$ instead of $\mathcal{E}(L)$ did not affect the form of the basis $\mathcal{B}_{c,1}$. More precisely, we have:
\begin{corintro*}[\textbf{\ref{animportantcor}}]
The set $\left\{a(z,t)\left|_{t=1}\,:\,\right. a(z,t)\in\mathcal{B}_{c,1}\right\}$ is linearly independent.
\end{corintro*}
Corollary \ref{animportantcor} allows us to transfer the graded nonlocal $\q$-vertex algebra structure from $V_{c,1}$
to $V_c=(V_{c,1})\left|_{t=1}\right. \subset\mathcal{E}(L)$. Hence we obtain the construction of the 
 graded nonlocal $\q$-vertex algebra $V_c$ generated by $x(z)\in\mathcal{E}(L)$.

In the end, we show that, for a suitably defined character $\ch_{\q}$, we have
\begin{thm*}[\textbf{\ref{charformula}}]
\begin{equation*}
\textstyle\ch_{\q} V_{c,1} =\displaystyle\sum_{r\geq 0} \frac{q^{r^2}}{(1-q)(1-q^2)\cdots(1-q^r)}c^r.
\end{equation*}
\end{thm*}
It is not surprising that for $c>1$ the character formula does not coincide with the classical case.
The analogous basis for $\widehat{\mathfrak{sl}}_2$  is derived using integrability \eqref{integrabilityrelation}
and commutativity of vertex operator products, while  vertex operator products for $U_q(\widehat{\mathfrak{sl}}_2)$ are no longer commutative, so only quantum
integrability \eqref{dm7} can be used.
Finally, in view of the close connection between the representation theory of  affine Kac-Moody Lie algebras and the Rogers-Ramanujan-type identities,
it is worth noting that the character formula from Theorem \ref{charformula} equals the left-hand side in 
\begin{align*}
&\sum_{r\geq 0} \frac{q^{r^2}}{(1-q)(1-q^2)\cdots(1-q^r)}c^r\\
&\quad=\left(
1+\sum_{s\geq 0} (-1)^s (1-cq^{2s}) c^{2s} q^{s(5s-1)/2} \frac{(1-cq)\cdots (1-cq^{s-1})}{(1-q)\cdots (1-q^s)}
\right)\prod_{r\geq 1}\frac{1}{1-c q^r},
\end{align*}
one of the standard identities which can be used to derive both Rogers-Ramanujan identities (\cite{Hardy}, cf. also \cite{Andrews}).

 \section{Preliminaries}\label{preliminaries}

 \subsection{Quantum calculus}
This subsection contains some elementary notions of quantum calculus.
For more details the reader may consult \cite{qKac}. 
 Fix an indeterminate $\q$. 
 For any two integers $m$ and
 $l$, $l\geq 0$, define $\q$-integers, $\q$-factorials, $\q$-binomial coefficients:
 \begin{align}
 &[m]_{\q}=\frac{\q^m -1}{\q-1}=1+\q+...+\q^{m-1}; \label{q61}\\
 &[0]_{\q} !=1,\quad [l+1]_{\q}!=[l+1]_{\q}[l]_{\q}\cdots[1]_{\q};\label{q62}\\
 &\gauss{m}{l}_{\q}=\frac{[m]_{\q}[m-1]_{\q}\cdots [m-l+1]_{\q}}{[l]_{\q}!}.\label{q63}
 \end{align}
 Denote by $z_0$ and $z$  two
 noncommutative variables satisfying 
 \begin{equation}\label{wz}
 z_0 z=\q zz_0.
 \end{equation}
 In the rest of this paper we shall assume that all  formal variables are commutative, unless stated otherwise (as above).
 We have the following $\q$-analogue of the binomial theorem:
  \begin{pro}\label{expand}
  For every integer $m$ and variables $z_0 ,z$ satisfying \eqref{wz} we have
  \begin{equation*}
 (z + z_0 )^m=\sum_{l\geq 0} \gauss{m}{l}_{\q} z^{m-l}z_{0}^{l}.
 \end{equation*}
  \end{pro}
  
  Let $V$ be a vector space over the field $\mathbb{C}(\q)$ and let
  $a(z)\in V[[z^{\pm 1}]]$ be an arbitrary Laurent series. Define
  {\em$\q$-derivation} of $a(z)$ as
  \begin{equation*}
  \frac{d_{\q}}{d_{\q} z}a(z)=\frac{a(z\q)-a(z)}{z(\q-1)}\in
  V [[z^{\pm 1}]].
  \end{equation*}
 In order to simplify our notation we will denote the $n$th
 $\q$-derivation of $a(z)$ as $a^{[n]}(z)$.
 The operator $\frac{d_{\q}}{d_{\q}z}$ is obviously a linear operator and it satisfies the general Leibniz rule:
 
 \begin{pro}\label{qLeibniz}
 For every nonnegative integer $m$ and $a(z),b(z)\in\om(V,V((z)))$  
 satisfying $a(z_1)b(z)\in\om(V,V((z_1,z)))$
 we have
 \begin{equation*}
 \left(a(z)b(z)\right)^{[m]}=\sum_{l=0}^{m}\gauss{m}{l}_{\q}a^{[l]}(z)b^{[m-l]}(z\q^l).
 \end{equation*}
 \end{pro}

\subsection{Quantum affine algebra \texorpdfstring{$U_{q}(\widehat{\mathfrak{g}})$}{Uq(g)}} \label{pre:1:2}
 First, we recall some facts from the theory of affine Kac-Moody Lie algebras (see \cite{Kac} for  details). 
  Let $\widehat{A}=(a_{ij})_{i,j=0}^{n}$ be a generalized Cartan matrix of affine type and let $S=\diag(s_{0},s_1,\ldots,s_n)$ be a diagonal matrix
  of relatively prime positive integers such that the matrix $S\widehat{A}$ is symmetric. 
 Let $\widehat{\mathfrak{t}}$ be a vector space over $\mathbb{C}(q^{1/2})$
 with a basis $\left\{\alpha^{\vee}_{0},\alpha^{\vee}_{1},\ldots,\alpha^{\vee}_{n},d\right\}$. Denote by $\alpha_{0}$, $\alpha_{1}$, \ldots, $\alpha_{n}$ linear functionals from $\widehat{\mathfrak{t}}^{*}$  such that
 $\alpha_{i}(\alpha_{j}^{\vee})=a_{ji}$ and $\alpha_{i}(d)=\delta_{i0}$
 for $i,j=0,1\ldots,n$.
 Define the set of simple roots, $\widehat{\Pi}=\left\{\alpha_{0},\alpha_1,...,\alpha_n\right\}$ and the set of simple coroots,  $\widehat{\Pi}^{\vee}=\left\{\alpha_{0}^{\vee},\alpha_1^{\vee},...,\alpha_n^{\vee}\right\}$.
 Denote by 
 $\widehat{\mathfrak{g}}$
 the affine Kac-Moody Lie algebra associated with the matrix $\widehat{A}$.
Let $\Lambda_{0}, \Lambda_{1}, \ldots, \Lambda_{n}$ be the fundamental weights, i.e. the elements of $\widehat{\mathfrak{t}}^{*}$ satisfying
 $$\Lambda_{i}(\alpha_{j}^{\vee})=\delta_{ij},\quad \Lambda_{i}(d)=0,
\quad i,j=0,1,\ldots,n.$$ 
Imaginary roots of $\widehat{\mathfrak{g}}$ are integer multiples of
 $\delta=d_{0}\alpha_{0}+d_{1}\alpha_{1}+\ldots+d_{n}\alpha_{n}\in \widehat{\mathfrak{t}}^{*},$
 where integers  $d_i$ are given in \cite{Kac}.
Define the integral dominant weight as any nonzero element $\Lambda$ of  the free Abelian group generated by $\Lambda_{0}, \Lambda_{1}, \ldots, \Lambda_{n},\delta/d_0$, which satisfies $\Lambda(\alpha_{i}^{\vee})\geq 0$ for $i=0,1,...,n$.
 The invariant symmetric bilinear form on $\widehat{\mathfrak{t}}^{*}$ is given by
 $(\alpha_{i},\alpha_{j})=s_{i}a_{ij}$, $(\delta,\alpha_{i})=(\delta,\delta)=0$, $i,j=0,1,\ldots,n$.
 
 Denote by $\mathfrak{g}$ a simple Lie algebra associated with the Cartan matrix 
 $A=(a_{ij})_{i,j=1}^{n}$.  Let $\mathfrak{t}\subset\widehat{\mathfrak{t}}$ be a Cartan subalgebra of  
 $\mathfrak{g}$, which is generated by the elements $\alpha_{1}^{\vee}$, $\alpha_{2}^{\vee}$, \ldots, $\alpha_{n}^{\vee}$.
Let
 $Q=\bigoplus_{i=1}^{n}\mathbb{Z}\alpha_{i}\subset\mathfrak{t}$ be the classical root lattice and  $P=\bigoplus_{i=1}^{n}\mathbb{Z}\lambda_{i}\subset\mathfrak{t}^{*}$ the classical weight lattice, where elements $\lambda_{i}\in\mathfrak{t}^{*}$ 
 satisfy $\lambda_{i}(\alpha_{j}^{\vee})=\delta_{ij}$ for $i,j=1,2,\ldots,n$.
 
 In this paper, we will mostly use $\q$-numbers defined in \eqref{q61}--\eqref{q63}.  However,
 the definition of quantum affine algebra is usually given in terms of (differently defined) $q$-numbers, so we recall  them as well.
 Fix an indeterminate $q$.
 For any two integers $m$ and $l$, $l\geq 0$,   define $q$-integers, $q$-factorials and $q$-binomial coefficients:
  \begin{align*}
  &[m]_q=\frac{q^{m}-q^{-m}}{q-q^{-1}};\\
  &[0]_q !=1,\quad [l+1]_q !=[l+1][l]\cdots[1];\\
& \gauss{m}{l}_q=\frac{[m][m-1]\cdots [m-l+1]}{[l]!}.
\end{align*}
  
 We  present the Drinfeld realization (see \cite{D}) 
 of the quantum affine algebra $U_{q}(\widehat{\mathfrak{g}})$.
 \begin{defn}\label{drinfeld}
 Let $\widehat{\mathfrak{g}}$ be an untwisted affine Kac-Moody Lie algebra with Cartan matrix $\widehat{A}=(a_{ij})_{i,j=0}^{n}$.
 The quantum affine algebra $U_{q}(\widehat{\mathfrak{g}})$ is the associative algebra over $\mathbb{C}(q^{1/2})$ with unit $1$ generated by the 
 elements
 $x_{i}^{\pm}(k)$, $a_{i}(l)$, $K_{i}^{\pm 1}$, $\gamma^{\pm 1/2}$ and $q^{\pm d}$, $i=1,2,\ldots,n$, $k,l\in\mathbb{Z}$, $l\neq 0$, subject to the following relations:
 {\allowdisplaybreaks
 \begin{align}
 & [\gamma^{\pm 1/2},u]=0\textrm{ for all }u\in U_{q}(\widehat{\mathfrak{g}})_0,\tag{d1}\label{D1}\\
 & K_i K_j=K_j K_i,\quad K_i K_{i}^{-1}=K_{i}^{-1}K_{i}=1,\tag{d2}\label{D2}\\
 & [a_{i}(k),a_{j}(l)]=\delta_{k+l\hspace{2pt}0}\frac{[a_{ij}k]_{q_{i}}}{k}\frac{\gamma^{k}-\gamma^{-k}}{q_{j}-q_{j}^{-1}},\tag{d3}\label{D3}\\
 & [a_{i}(k),K_{j}^{\pm 1}]=[q^{\pm d},K_{j}^{\pm 1}]=0,\tag{d4}\label{D4}\\
 & q^{d}x_{i}^{\pm}(k)q^{-d}=q^{k}x_{i}^{\pm}(k),\quad q^{d}a_{i}(l)q^{-d}=q^{k}a_{i}(l),\tag{d5}\label{D5}\\
 & K_{i}x_{j}^{\pm}(k)K_{i}^{-1}=q^{\pm (\alpha_{i},\alpha_{j})}x_{j}^{\pm }(k) ,\tag{d6}\label{D6}\\
 & [a_{i}(k),x_{j}^{\pm}(l)]=\pm\frac{[a_{ij}k]_{q_{i}}}{k}\gamma^{\mp |k|/2}x_{j}^{\pm}(k+l),\tag{d7}\label{D7}\\
 & x_{i}^{\pm}(k+1)x_{j}^{\pm}(l)-q^{\pm(\alpha_i,\alpha_j)}x_{j}^{\pm}(l)x_{i}^{\pm}(k+1)\nonumber\\
 &\hspace{20pt}=q^{\pm(\alpha_i,\alpha_j)}x_{i}^{\pm}(k)x_{j}^{\pm}(l+1)-x_{j}^{\pm}(l+1)x_{i}^{\pm}(k),\tag{d8}\label{D8}\\
 & [x_{i}^{+}(k),x_{j}^{-}(l)]=\frac{\delta_{ij}}{q_{i}-q^{-1}_{i}}\left(\gamma^{\frac{k-l}{2}}\psi_{i}(k+l)-\gamma^{\frac{l-k}{2}}\phi_{i}(k+l)\right),\tag{d9}\label{D9}\\
 & \sym_{l_1,l_2,\ldots,l_m}\sum_{s=0}^{m}(-1)^{s}\gauss{m}{s}_{q_i}x_{i}^{\pm}(l_{1})\cdots x_{i}^{\pm}(l_{s})x_{j}^{\pm}(k)x_{i}^{\pm}(l_{s+1})\cdots x_{i}^{\pm}(l_{m})=0, \quad\textrm{for }i\neq j,\tag{d10}\label{Da}
 \end{align}}
 where $m=1-a_{ij}$, $q_i =q^{s_i}$ and the elements $\phi_{i}(-r)$ and $\psi_{i}(r)$, $r\in\mathbb{Z}_{\geq 0}$, are given by \\
 \begin{align*}
 & \phi_{i}(z)=\sum_{r=0}^{\infty}\phi_{i}(-r)z^{r}=K_{i}^{-1}\exp\left(-(q_{i}-q_{i}^{-1})\sum_{r=1}^{\infty}a_{i}(-r)z^{r}\right),\\
 & \psi_{i}(z)=\sum_{r=0}^{\infty}\psi_{i}(r)z^{-r}=K_{i}\exp\left((q_{i}-q_{i}^{-1})\sum_{r=1}^{\infty}a_{i}(r)z^{-r}\right).
 \end{align*}
 \end{defn}
 If $q_{i}=q$ we will usually omit the index $q_i$ and write $[m]$ instead of $[m]_{q_{i}}$. Denote by $x_{i}^{\pm}(z)$ the series
 \begin{equation}\label{101_exp:series}
 x_{i}^{\pm}(z)=\sum_{r\in\mathbb{Z}}x_{i}^{\pm}(r)z^{-r-1}\in U_{q}(\widehat{\mathfrak{g}})[[z,z^{- 1}]].
 \end{equation}
 We shall continue to use the notation $x_{i}^{\pm}(z)$ for the action of  (\ref{101_exp:series}) on an arbitrary
 $U_{q}(\widehat{\mathfrak{g}})$-module $V$:
 $$x_{i}^{\pm}(z)=\sum_{r\in\mathbb{Z}}x_{i}^{\pm}(r)z^{-r-1}\in (\ndo V)[[z,z^{- 1}]].$$
 
 Drinfeld gave the following Hopf algebra structure for his realization of   $U_{q}(\widehat{\mathfrak{sl}}_{n+1})$  (cf. also \cite{DI}):
{\allowdisplaybreaks
 \begin{align*}
 &\Delta(q^{c/2})=q^{c/2}\otimes q^{c/2},\\
 &\Delta(x_{i}^{+}(z))=x_{i}^{+}(z)\otimes 1 + \phi_{i}(zq^{c_{1}/2})\otimes x_{i}^{+}(zq^{c_{1}}),\\
 &\Delta(x_{i}^{-}(z))=1\otimes x_{i}^{-}(z)+x_{i}^{-}(zq^{c_{2}})\otimes\psi_{i}(zq^{c_{2}/2}),\\
 &\Delta(\phi_{i}(z))=\phi_{i}(zq^{-c_{2}/2})\otimes \phi_{i}(zq^{c_{1}/2}),\\
 &\Delta(\psi_{i}(z))=\psi_{i}(zq^{c_{2}/2})\otimes \psi_{i}(zq^{-c_{1}/2}),\\
 &\varepsilon(q^c)=1, \quad \varepsilon(x^{\pm}_{i}(z))=0,\quad\varepsilon(\phi_{i}(z))=\varepsilon(\psi_{i}(z))=1;\\
 & S(q^{c})=q^{-c},\\
 & S(x_{i}^{+}(z))=-\phi_{i}(zq^{-c/2})^{-1}x_{i}^{+}(zq^{-c}),\\
 & S(x_{i}^{-}(z))=-x_{i}^{-}(zq^{-c})\psi_{i}(zq^{-c/2})^{-1},\\
 & S(\phi_{i}(z))=\phi_{i}(z)^{-1},\quad S(\psi_{i}(z))=\psi_{i}(z)^{-1},
 \end{align*}}where $i=1,2,...,n$ and $q^{c_{1}}$ denotes the action of the center $\gamma=q^{c}$ on the first tensor factor while $q^{c_{2}}$ the action  on the second tensor factor.
 For $l\geq 1$ set
 \begin{equation*}
 \Delta^{(0)}=1\quad\textrm{ and }\quad\Delta^{(l)}=(\underbrace{1\otimes\cdots\otimes 1}_{\mbox{$l-1$}} \otimes \Delta)\Delta^{(l-1)}.
 \end{equation*}
 The coproduct formula applied on the tensor product of $c$ integrable highest weight modules gives
 \begin{equation}\label{copro1}
 \Delta^{(c-1)}(x_{i}^{+}(z))=\sum_{l=1}^{c}x_{i}^{+(l)}(z),
 \end{equation}
 where
 \begin{equation}\label{copro2}
 x_{i}^{+(l)}(z)=\underbrace{\phi_{i}(zq^{\frac{1}{2}})\otimes\phi_{i}(zq^{\frac{3}{2}})\otimes\cdots\otimes\phi_{i}(zq^{l-\frac{3}{2}})}_{\mbox{$l-1$}}\otimes x_{i}^{+}(zq^{l-1})\otimes\underbrace{1\otimes\cdots\otimes 1}_{\mbox{$c-l$}}.
 \end{equation}

 \subsection{Representations of \texorpdfstring{$U_{q}(\widehat{\mathfrak{sl}}_{2})$}{Uq(sl2)}}

 The algebra $U_{q}(\widehat{\mathfrak{sl}}_{2})$ is generated by the elements
  $x_{1}^{\pm}(k)$, $a_{1}(l)$, $K_{1}^{\pm 1}$, $\gamma^{\pm 1/2}$ and $q^{\pm d}$, $k,l\in\mathbb{Z}$, $l\neq 0$.
 In order to simplify the notation, we will omit  index "$1$" and write $x^{\pm}(k)$, $a(l)$, $K^{\pm 1}$, $\phi(z)$
 instead of $x_{1}^{\pm}(k)$, $a_{1}(l)$, $K_{1}^{\pm 1}$, $\phi_1 (z)$.
  We recall  Frenkel-Jing realization of  the integrable highest weight
   $U_{q}(\widehat{\mathfrak{sl}}_{2})$-modules $L(\Lambda_{0})$ and $L(\Lambda_1)$ (see \cite{FJ}).
 The Heisenberg algebra  $U_{q}(\widehat{\mathfrak{h}})$ of level $1$ is generated by the elements $a(l)$,  $l\in\mathbb{Z}\setminus\left\{0\right\}$, and the central element $\gamma^{\pm 1}=q^{\pm c}$ subject to the relations 
 \begin{equation}\label{104_heisenberg}
 [a(r),a(s)]=\delta_{r+s\hspace{2pt}0}\frac{[2r][r]}{r},\quad
  r,s\in\mathbb{Z}\setminus\left\{0\right\}.
 \end{equation}
 Algebra $U_{q}(\widehat{\mathfrak{h}})$ has a natural realization on the space $\sym(\widehat{\mathfrak{h}}^{-})$
 of the symmetric algebra generated by the elements $a(-r)$, 
   $r\in\mathbb{Z}_{>0}$, via the following rule:
   \begin{align*}
 \gamma^{\pm 1}\hspace{5pt}&\ldots\hspace{5pt}\textrm{multiplication by }q^{\pm 1},\\
 a(r)\hspace{5pt}&\ldots\hspace{5pt}\textrm{differentiation operator subject to (\ref{104_heisenberg})},\\
 a(-r)\hspace{5pt}&\ldots\hspace{5pt}\textrm{multiplication by the element }a(-r).
 \end{align*}
 Denote the resulted level $1$ irreducible $U_{q}(\widehat{\mathfrak{h}})$-module 
  by $M(1)$. Define the following operators on $M(1)$:
  \begin{align*}
  &E_{-}^{\pm}(a,z)=\exp\left(\mp\sum_{r\geq 1}\frac{q^{\mp r/2}}{[r]}a(-r)z^{r}\right),
   \quad E_{+}^{\pm}(a,z)=\exp\left(\pm\sum_{r\geq 1}\frac{q^{\mp r/2}}{[r]}a(r)z^{-r}\right).
 \end{align*}
 
 Let $\mathbb{C}\left\{Q\right\}$  be the  group algebra of the classical weight lattice $Q=\mathbb{Z}\alpha_1$ generated by $e^{\alpha}$, $\alpha\in Q$. The space
 $\mathbb{C}\left\{P\right\}=\mathbb{C}\left\{Q\right\} \oplus \mathbb{C}\left\{Q\right\}e^{\lambda_1}$
 is a $\mathbb{C}\left\{Q\right\}$-module.
  For $\alpha\in Q$ define an action  $z^{\partial_\alpha}$ on $\mathbb{C}\left\{P\right\}$ by
  $z^{\partial_\alpha}e^{\beta}=z^{(\alpha,\beta)}e^{\beta}.$
  Set 
 \begin{equation*}
 L_0 =M(1)\otimes\mathbb{C}\left\{Q\right\},\qquad L_1 =M(1)\otimes\mathbb{C}\left\{Q\right\}e^{\lambda_1}.
 \end{equation*}

 \begin{thm}[\cite{FJ}]\label{frenkeljing}
  By the action 
 \begin{align*}
 x^{\pm}(z)&=E_{-}^{\pm}(-a,z)E_{+}^{\pm}(-a,z)\otimes e^{\pm\alpha_1}z^{\pm\partial_{\alpha_1}},
 \end{align*}
 the space $L_i$, $i=0,1$, becomes
 the integrable highest weight $U_{q}(\widehat{\mathfrak{sl}}_{2})$-module    with the highest weight $\Lambda_i$.
 \end{thm}
 The following theorem is a special case of the, so-called quantum integrability, which was found by  Ding and  Miwa.
 \begin{thm} [\cite{DM}]
 On every  level $c$ integrable $U_{q}(\widehat{\mathfrak{sl}}_{2})$-module
 \begin{equation}\label{ding_miwa}
 x^{\pm}(z_{1})x^{\pm}(z_{2})\cdots x^{\pm}(z_{c+1})=0\quad\text{if}\quad z_{1}/z_{2}=z_{2}/z_{3}=\ldots=z_{c}/z_{c+1}=q^{\mp 2}.
 \end{equation}
\end{thm}

In the end, following
 \cite{Koyama}
we define the   vertex operator $\mathcal{Y}(z)$ on the space
 $V=M(1)\otimes \mathbb{C}\left\{P\right\}$:
 \begin{align*}
 &\mathcal{E}_{-}(z)=\exp\left(\sum_{r\geq 1}\frac{q^{r/2}}{[2r]}a(-r)z^{r}\right),\quad
 \mathcal{E}_{+}(z)=\exp\left(-\sum_{r\geq 1}\frac{q^{r/2}}{[2r]}a(r)z^{-r}\right);\\
 &\mathcal{Y}(z)=\mathcal{E}_{-}(z)\mathcal{E}_{+}(z)
 \otimes e^{\lambda_1}(-z)^{\partial_{\lambda_1}}\in\om(V,V((z^{1/2}))).
 \end{align*}
 The relations in the next proposition can be verified by a direct calculation. 
 \begin{pro}\label{relations}
 The following relations hold on $V$:
 {\allowdisplaybreaks\begin{align}
 & x^{+}(z_1)x^{+}(z_2)=(z_1-z_2)(z_1-q^{-2}z_2):x^{+}(z_1)x^{+}(z_2):,\label{r1}\\
  &x^{+}(z_{1})\phi(z_{2}q^{1/2})=q^{2}\frac{1-q^{-2}\frac{z_{2}}{z_{1}}}{1-q^{2}\frac{z_{2}}{z_{1}}}\phi(z_{2}q^{1/2})x^{+}(z_{1}),\label{normal1}\\
 & x^{+}(z_1)\mathcal{E}_{-}(z_2)=\left(1-\frac{z_2}{z_1}\right)\mathcal{E}_{-}(z_2)x^{+}(z_1),\label{komutiranjesE-}\\
  & \phi (z_1)\mathcal{E}_{-}(z_2)=\mathcal{E}_{-}(z_2)\phi (z_1).\label{785}
 \end{align}}
 \end{pro}

\section{Graded nonlocal \texorpdfstring{$\q$}{q}-vertex algebras}\label{gkva}

\subsection{Definition}
Let $V$ be an arbitrary  vector space over the field $\mathbb{F}\supseteq \mathbb{C}(\q)$ of characteristic zero 
 and $1_V$ the identity $V\to V$.  

\begin{defn}\label{kva}
A {\em graded nonlocal $\q$-vertex algebra} is a triple $(V,Y,1)$, where  
\begin{equation}\tag{v0}\label{VVV}
V= \coprod_{u\in\mathbb{Z}} V_{(u)}
\end{equation}
 is a
$\mathbb{Z}$-graded vector space
 equipped with a
 linear map
\begin{align}
Y(\cdot,z_0)\colon V &\to (\ndo V)[[z^{\pm 1}_{0}]]\nonumber\\
a&\mapsto Y(a,z_0 )=\sum_{r\in\mathbb{Z}}a_{r}z^{-r-1}_{0}\tag{v1}\label{V0}
\end{align}
\noindent and  with a distinguished vector $\vac\in V$ such that the
following conditions hold: For every $a,b,c\in V$ and $s,u\in\mathbb{Z}$
\begin{align}
&a_r b=0\quad\text{for  sufficiently large integer }r;\tag{v2}\label{V1}\\
&Y(\vac,z_0)=1_{V};\tag{v3}\label{V2}\\
&Y(a,z_0)\vac\in V[[z_0]]\quad\text{and}\quad \lim_{z_0\to
0}Y(a,z_0)\vac=a;\tag{v4}\label{V3}\\
&a_r V_{(s)}\subseteq V_{(s+u-r-1)}\quad\text{if }a\in V_{(u)};\tag{v5}\label{V4}\\
& Y(a,z_0 +z_2)Y(b,z_2 \q^{u}) c=Y(Y(a,z_0)b,z_2)c\quad\text{if }a(z)\in V_{(u)}\text{ and}\tag{v6}\label{V5}\\
& z_2 z_0=\q z_0 z_2.\tag{v7}\label{commcons}
\end{align}
\end{defn}
The definition requires some further explanations.
We use the following conventions in  \eqref{V5}. On the left-hand side 
we assume that $z_0 +z_2$ appears to the right of $z_2$:
\begin{align}\label{cfleft}
 Y(a,z_0 +z_2)Y(b,z_2\q^u) c=\sum_{r\in\mathbb{Z}}\sum_{s\in\mathbb{Z}} a_r ( b_s c) (z_2\q^u)^{-s-1} (z_0 +z_2)^{-r-1}.
\end{align}
The expression $(z_0 +z_2)^{-r-1}$ is expanded in nonnegative powers of $z_2$ using Proposition \ref{expand}.
On  the right-hand side of \eqref{V5}, we  assume that the variable $z_0$ appears to the left of $z_2$:
\begin{align}
&Y(Y(a,z_0)b,z_2)c=Y(\sum_{r\in\mathbb{Z}}a_r bz_0^{-r-1},z_2)c=\sum_{s\in\mathbb{Z}} \sum_{r\in\mathbb{Z}}(a_r b)_s cz_0^{-r-1} z_2^{-s-1}.\label{cfright}
\end{align}
As a consequence of \eqref{V5}, using \eqref{cfleft} and \eqref{cfright}, we get
\begin{equation}\label{asocc}
(a_r b)_s c =\sum_{l\geq 0}\gauss{l-r-1}{l}_\q \q^{(s+l+1)(r-u+1)}a_{r-l}(b_{s+l}c)\quad\text{for all }r,s\in\mathbb{Z}.
\end{equation}

\subsection{Construction}

Let $L$ be an arbitrary  vector space over the field $\mathbb{F}\supseteq \mathbb{C}(\q)$ of characteristic zero.
We will use the following notation:
\begin{align*}
 \mathcal{E}(L)=\om(L,L((z))),\qquad \mathcal{E}(L)_t=\om(L,L((z)))\otimes \mathbb{F}[t].
 \end{align*}
Our goal is to construct graded nonlocal $\q$-vertex algebras generated by the subsets of $\mathcal{E}(L)_t$.
For  $A(z_1,...,z_n)\in \om(L,L((z_1,...,z_n)))$ we will denote
by $\lim_{z_j \to z} A(z_1,...,z_n)$ or $\lim_{z_1,...,z_n \to z} A(z_1,...,z_n)$ an element of the space
 $\om(L,L((z)))$ which is
 obtained from $A(z_1,...,z_n)$ by replacing all the variables  $z_1,...,z_n$ by $z$.
Notice that $\lim_{z_j \to z} B(z_1,...,z_n)$ need not  be well-defined for an arbitrary $B(z_1,...,z_n)\in \om(L,L[[z_1^{\pm 1},...,z_n^{\pm 1}]])$.
First, we recall a few well-known technical results which will be often used in this section.

\begin{lem}\label{mp} Let  $A(z_1,z)\in\om(L,L((z_1,z)))$, $B(z_2,z_1,z)\in\om(L,L((z_2,z_1,z)))$ and $p(z)\in\mathbb{F}(\q)[z]$. Then
\begin{enumerate}[(a)]
\item\label{mp2} $\displaystyle\lim_{z_1\to z}  A(z_1,z) =\lim_{z_1\to z} \lim_{z\to z_1} A(z_1,z)$;
\item\label{mp3} $\displaystyle\lim_{z_2,z_1\to z} B(z_2,z_1,z) = \lim_{z_2\to z}\lim_{z_1\to z}  B(z_2,z_1,z)$;
\item\label{mp1} $\displaystyle\lim_{z_1\to z} p(z/z_1) A(z_1,z) = \lim_{z_1\to z} p(z/z_1) \lim_{z_1\to z} A(z_1,z)$.
\end{enumerate}
\end{lem}

\begin{prf}
Equality \eqref{mp2} is trivial while equalities \eqref{mp3} and \eqref{mp1} follow from more general results proved in \cite[Section 5]{MP}.
\end{prf}

Denote by $\mathbb{F}(\q)[z]_0$ the following set of polynomials:
$$\mathbb{F}(\q)[z]_0 = \left\{p(z)\in\mathbb{F}(\q)[z]\, :\, p(1)=1\text{ and }p(\q^n)\neq 0 \text{ for all }  n\in\mathbb{Z}_{\geq 0} \right\}.$$
Obviously, $\mathbb{F}(\q)[z]_0$ is closed under multiplication.
We shall consider the space
$$\om(L,L((z_1,z_2,\ldots,z_m)))_t = \om(L,L((z_1,z_2,\ldots,z_m)))\otimes\mathbb{F}(\q)[t].$$
For a homogeneous element $a(z)\otimes t^{\alpha}\in \mathcal{E}(L)_t$ we  write 
\begin{equation}\label{newlabel}
\wt (a(z)\otimes t^{\alpha}) =\deg_t (a(z)\otimes t^{\alpha})=\alpha.
\end{equation}
The function $\wt$ will define a gradation on the graded nonlocal $\q$-vertex algebras which will be constructed in this section.

\begin{defn}\label{qkomm}
A sequence  $(a_1 (z)\otimes t^{\alpha_1},\ldots,a_m (z)\otimes t^{\alpha_m})$ in
$\mathcal{E}(L)_t$   is said to be {\em quasi-commutative} if there exists a polynomial $p(z)\in\mathbb{F}(\q)[z]_0$ such that 
\begin{align}
\left(\prod_{1\leq i < j\leq m} p(z_j/z_i) \right)a_1(z_1)a_2 (z_2 \q^{\alpha_1})a_3 (z_3 \q^{\alpha_1+\alpha_2})\cdots a_m(z_m\q^{\alpha_1+...+\alpha_{m-1}})&\label{qcomm}\\
\quad\in\om(L,L((z_1,z_2,\ldots,z_m))).&\nonumber
\end{align}
In general, a sequence $(\sum_{j=1}^{N_1}a_{1,j} (z)\otimes t^{\alpha_{1,j}},\ldots,\sum_{j=1}^{N_m}a_{m,j} (z)\otimes t^{\alpha_{m,j}})$ in
$\mathcal{E}(L)_t$ is said to be {\em quasi-commutative} if for every choice of $(j_1,...,j_m)$, where $1\leq j_k\leq N_k$, $k=1,...,m$, the sequence
$(a_{1,j_1} (z)\otimes t^{\alpha_{1,j_1}},\ldots,a_{m,j_m} (z)\otimes t^{\alpha_{m,j_m}})$ is quasi-commutative.
\end{defn}

Notice that if  $(a_1 (z)\otimes g_1(t),\ldots,a_m (z)\otimes g_m(t))$ is quasi-commutative for some polynomials $g_j(t)$, $j=1,...,m$, then for all integers $\alpha_1\leq \alpha_2\leq ...\leq \alpha_m$ the sequence  $(a_1 (z\q^{\alpha_1})\otimes g_1(t),\ldots,a_m (z\q^{\alpha_m})\otimes g_m(t))$ is also quasi-commutative.

\begin{defn}\label{produkti}
Let $(a(z,t),b(z,t))=(\sum_{i=1}^{M}a_i (z)\otimes t^{\alpha_i},\sum_{j=1}^{N}b_j(z)\otimes t^{\beta_j})$ be a quasi-commutative pair in $\mathcal{E}(L)_t$.
 For an integer
$r$ we define 
$$a(z,t)_{-r-1}b(z,t)\in (\ndo L)[[z^{\pm 1}]]_t=(\ndo L)[[z^{\pm 1}]]\otimes\mathbb{F}(\q)[t]$$
 by
\begin{equation}\label{rtiprodukti}
a(z,t)_{-r-1}b(z,t) = 
\sum_{i=1}^{M}\sum_{j=1}^{N}\left(a_i(z)\otimes t^{\alpha_i}\right)_{-r-1} \left(b_j(z)\otimes t^{\beta_j}\right)
\end{equation}
where, for $r\geq 0$
\begin{equation}\label{rtiprodukti2}
 (a_i(z)\otimes t^{\alpha_i} )_{-r-1} (b_j(z)\otimes t^{\beta_j})=\lim_{z_1\to z} \frac{p_{ij}(z/z_1)}{[r]_\q !}a_i^{[r]}(z_1)b_j(z\q^{\alpha_i +r})\otimes t^{\alpha_i +\beta_j  +r}
\end{equation}
and $p_{ij}(z)\in\mathbb{F}(\q)[z]_0$ is any polynomial satisfying
\begin{equation}\label{rtipolinom}
p_{ij}(z/z_1)a_i^{[r]}(z_1)b_j(z\q^{\alpha_i +r})\in \om(L,L((z_1,z))),
\end{equation}
while  for $r<0$ we set
$$a(z,t)_{-r-1}b(z,t) =0.$$
\end{defn}

Notice that the polynomial $p_{ij}(z)$ in Definition \ref{produkti} does depend on the choice of the integer $r$. At this point, it is not clear whether this definition, more precisely  \eqref{rtiprodukti2}, depends on the choice of the polynomial  $p_{ij}(z)$. This issue is resolved in the next proposition.

\begin{pro}\label{proposition_independence}
Definition \ref{produkti} is independent of the choice of the polynomials $p_{ij}(z)\in\mathbb{F}(\q)[z]_0$ satisfying \eqref{rtipolinom}.
\end{pro}

\begin{prf}
The proof goes similarly as the proof of \cite[Lemma 3.3]{Li0}.
It is sufficient to consider only homogeneous elements $a(z)\otimes t^{\alpha},b(z)\otimes t^{\beta}\in\mathcal{E}(L)_t$.
Let $p_k (z)\in\mathbb{F}(\q)[z]_0$, $k=1,2$, be two polynomials satisfying
\begin{equation*}
p_k (z/z_1)a^{[r]}(z_1)b(z\q^{\alpha +r})\in \om(L,L((z_1,z)))
\end{equation*}
for some $r\geq 0$. Using Lemma \ref{mp} we get
\allowdisplaybreaks
\begin{align*}
&\lim_{z_1\to z} \frac{p_1 (z/z_1)}{[r]_\q !}a^{[r]}(z_1)b(z\q^{\alpha +r})\otimes t^{\alpha+\beta +r}\\
=&p_2 (1) \lim_{z_1\to z} \frac{p_1 (z/z_1)}{[r]_\q !}a^{[r]}(z_1)b(z\q^{\alpha +r})\otimes t^{\alpha+\beta +r}\\
=& \lim_{z_1\to z}p_2 (z/z_1) \frac{p_1 (z/z_1)}{[r]_\q !}a^{[r]}(z_1)b(z\q^{\alpha +r})\otimes t^{\alpha+\beta +r}\\
=& \lim_{z_1\to z}p_1 (z/z_1) \frac{p_2 (z/z_1)}{[r]_\q !}a^{[r]}(z_1)b(z\q^{\alpha +r})\otimes t^{\alpha+\beta +r}\\
=&p_1 (1) \lim_{z_1\to z} \frac{p_2 (z/z_1)}{[r]_\q !}a^{[r]}(z_1)b(z\q^{\alpha +r})\otimes t^{\alpha+\beta +r}\\
=&\lim_{z_1\to z} \frac{p_2 (z/z_1)}{[r]_\q !}a^{[r]}(z_1)b(z\q^{\alpha +r})\otimes t^{\alpha+k +r},
\end{align*}
as required.
\end{prf}

The following generalization of Proposition \ref{proposition_independence} can be proved analogously.

\begin{lem}
Let $(a_1(z)\otimes t^{\alpha_1}, ... ,a_m(z)\otimes t^{\alpha_m})$ be a quasi-commutative sequence in $\mathcal{E}(L)_t$ and let $p_k (z)\in\mathbb{F}(\q)[z]_0$, $k=1,2$,  be polynomials satisfying 
\begin{align*}
A_k(z_1,...,z_m):=\left(\prod_{1\leq i < j\leq m} p_k(z_j/z_i) \right)a_1(z_1)a_2 (z_2 \q^{\alpha_1})a_3 (z_3 \q^{\alpha_1+\alpha_2})\cdots a_m(z_m\q^{\alpha_1+...+\alpha_{m-1}})&\\
\in\om(L,L((z_1,z_2,\ldots,z_m)))\quad\text{for }k=1,2.&
\end{align*}
Then
$$\lim_{z_1,...,z_m\to z}A_1(z_1,...,z_m)=\lim_{z_1,...,z_m\to z}A_2(z_1,...,z_m).$$
\end{lem}

Let $(a(z)\otimes t^{\alpha},b(z)\otimes t^{\beta})$ be a quasi-commutative pair in $\mathcal{E}(L)_t$. Set
\begin{equation}\label{sydney}
Y(a(z)\otimes t^{\alpha},z_0)(b(z)\otimes t^{\beta}) =\sum_{r\in\mathbb{Z}}(a(z)\otimes t^{\alpha})_{-r-1}(b(z)\otimes t^{\beta}) z_0^r
\end{equation}
and then extend $Y$ by linearity. Set 
\begin{equation}\label{vacuum}
\vac =1_{\mathcal{E}(L)}\otimes 1\in\mathcal{E}(L)_t .
\end{equation}
The following properties of the $r$th products follow directly from \eqref{newlabel}, \eqref{sydney}, \eqref{vacuum} and Definition \ref{produkti}.

\begin{kor}\label{lotsofthings}
Let $(a(z,t),b(z,t))$ and $(c(z)\otimes t^{\alpha},d(z)\otimes t^{\beta})$ be two quasi-commutative pairs in $\mathcal{E}(L)_t$
and $r$ an arbitrary integer.
\begin{enumerate}[(a)]
\item $a(z,t)_{r}b(z,t)\in\mathcal{E}(L)_t$  for $r\in\mathbb{Z}$;
\item The pairs $(a(z,t),\vac)$  and  $(\vac,a(z,t))$  are quasi-commutative;
\item $Y(\vac,z_0)=1_{\mathcal{E}(L)_t}$;
\item $Y(a(z,t),z_0)\vac\in \mathcal{E}(L)_t [[z_0]]$ and $ \lim_{z_0\to
0}Y(a(z,t),z_0)\vac=a(z,t)$;
\item If $e(z,t):=(c(z)\otimes t^{\alpha})_{r}(d(z)\otimes t^{\beta})\neq 0$, then  $\wt(e(z,t)) =\alpha+\beta -r-1$.
\end{enumerate}
\end{kor}

The following lemma is one of the key results in this section.

\begin{lem}\label{dong}
Let  $(a_1(z,t), ... ,a_m(z,t))$ and $(a_{k}(z,t),a_{k+1}(z,t))$ for some  $k=1,2,...,m-1$ be two quasi-commutative sequences in $\mathcal{E}(L)_t$. Then for every integer $r$    the sequence
\begin{align*}
&(a_1(z,t),...,a_{k-1}(z,t) ,a_{k}(z,t)_{-r-1} a_{k+1}(z,t)  ,
a_{k+2}(z,t),...  ,a_{m}(z,t) )
\end{align*}
 is quasi-commutative.
\end{lem}

\begin{prf}
Let $r\geq 0$.
It is sufficient to consider only the sequence of homogeneous elements
$(a_1(z)\otimes t^{\alpha_1}, ... ,a_m(z)\otimes t^{\alpha_m}).$
Let $p(z)\in\mathbb{F}(\q)[z]_0$ be a polynomial satisfying
\begin{align}
&p(z_{k+1}/z_k)a_{k}(z_k)a_{k+1}(z_{k+1}\q^{\alpha_k})\in\om(L,L((z_k, z_{k+1}))) ;\label{dongv}\\
&p(z_{k+1}/z_k)a_{k}(z_k\q^l)a_{k+1}(z_{k+1}\q^{r+\alpha_k})\in\om(L,L((z_k, z_{k+1})))\quad\text{for }l=0,1,...,r; \label{dong2}\\
&\left(\prod_{1\leq i< j\leq m} p(z_j / z_i)\right) a_1(z_1)a_2 (z_2 \q^{\alpha_1})\cdots a_m(z_m\q^{\alpha_1+...+\alpha_{m-1}})\in\om(L,L((z_1,..., z_{m}))).\label{dong4}
\end{align}
Note that \eqref{dong2} implies
$$p(z_{k+1}/z_k)a_{k}^{[r]}(z_k)a_{k+1}(z_{k+1}\q^{r+\alpha_k})\in\om(L,L((z_k, z_{k+1}))).$$
Since $(a_{k}(z)\otimes t^{\alpha_k})_{-r-1} (a_{k+1}(z) \otimes t^{\alpha_{k+1}})$ is a $\mathbb{F}(\q)[z^{- 1}]$-linear combination of 
$$(a_{k}(z\q^l)\otimes t^{\alpha_k +r})_{-1} (a_{k+1}(z) \otimes t^{\alpha_{k+1}}),\quad l=0,1,...,r,$$ 
it is sufficient to prove that the sequence
\begin{align}
&(a_1(z)\otimes t^{\alpha_1},...,a_{k-1}(z)\otimes t^{\alpha_{k-1}}  ,(a_{k}(z\q^l)\otimes t^{\alpha_k +r})_{-1} (a_{k+1}(z) \otimes t^{\alpha_{k+1}}),\nonumber\\
& \quad a_{k+2}(z)\otimes t^{\alpha_{k+2}} ,...  ,a_{m}(z)\otimes t^{\alpha_{m}}  )\label{dong8}
\end{align}
 is quasi-commutative for every $l=0,1,...,r$. First, \eqref{dong4} implies
\begin{equation}\label{dong6}
\left(\prod_{1\leq i< j\leq m} p(w_j / w_i)\right) a_1(w_1)a_2(w_2\q^{\alpha_1})\cdots a_m(w_m\q^{\alpha_1+...+\alpha_{m-1}})\in\om(L,L((z_1,..., z_{m}))),
\end{equation}
where $w_u = z_u$ for $u< k$, $w_k =z_k \q^{l}$ and $w_u = z_u \q^{r}$ for $u>k$.  
Notice that $p(w_j / w_i)$ is an element of $\mathbb{F}(\q)[z_j /z_i]_0$ for all $1\leq i< j\leq m$.
Expression
\eqref{dong6} can be written as
\begin{align}
&\left(\prod_{\substack{1\leq i< j\leq m\\ (i,j)\neq (k,k+1)}} p(w_j / w_i)\right) a_1(w_1)
\cdots a_{k-1}(w_{k-1}\q^{\alpha_1+...+\alpha_{k-2}})\nonumber\\
&\quad\cdot\left(p(w_{k+1} / w_k)a_k(w_k\q^{\alpha_1+...+\alpha_{k-1}})a_{k+1}(w_{k+1}\q^{\alpha_1+...+\alpha_{k}})\right)\nonumber\\
&\quad\quad\cdot a_{k+2}(w_{k+2}\q^{\alpha_1+...+\alpha_{k+1}})\cdots
 a_m(w_m\q^{\alpha_1+...+\alpha_{m-1}}).\label{dong7}
\end{align}
Applying \eqref{dong7}  on an arbitrary vector $v\in L$, we get a finite number of negative powers of $z_j$, $j=1,2,...,m$. Therefore,  applying the limit $\lim_{z_k \to z_{k+1}}$ on \eqref{dong7} and using \eqref{dongv} we get, up to a nonzero multiplicative scalar,
\begin{align}
&\left(\prod_{\substack{1\leq i< j\leq m\\ i,j\neq k}} p(w_j / w_i)\right)
\left(\prod_{j>k+1} p(w_j /(w_{k+1}\q^{l-r}) ) \right)
\left(\prod_{i<k}p((w_{k+1}\q^{l-r}) / w_i)\right)\label{dong_polinomial}\\
&\quad\cdot a_1(w_1)\cdots a_{k-1}(w_{k-1}\q^{\alpha_1+...+\alpha_{k-2}})\nonumber\\
&\quad\quad\cdot a_{k}(w_{k+1}\q^{\alpha_1+...+\alpha_{k-1}+l-r})_{-1} a_{k+1}(w_{k+1}q^{\alpha_1+...+\alpha_{k}})\label{dong_polinomial2} \\
&\quad\quad\quad\cdot a_{k+2}(w_{k+2}\q^{\alpha_1+...+\alpha_{k+1}})\cdots
 a_m(w_m\q^{\alpha_1+...+\alpha_{m-1}})\nonumber\\
&\quad\quad\quad\quad\in\om(L,L((z_1,...,z_{k-1},z_{k+1},... z_{m}))),\nonumber
\end{align}
where
\begin{align*}
& a_{k}(w_{k+1}\q^{\alpha_1+...+\alpha_{k-1}+l-r})_{-1} a_{k+1}(w_{k+1}q^{\alpha_1+...+\alpha_{k}})\\
&\quad =
a_{k}(z_{k+1}\q^{\alpha_1+...+\alpha_{k-1}+l})_{-1} a_{k+1}(z_{k+1}q^{\alpha_1+...+\alpha_{k}+r})
\end{align*}
in \eqref{dong_polinomial2} denotes the first tensor factor of
\begin{align*}
&(a_{k}(z_{k+1}\q^{\alpha_1+...+\alpha_{k-1}+l})\otimes t^{\alpha_k +r})_{-1} (a_{k+1}(z_{k+1}q^{\alpha_1+...+\alpha_{k-1}})\otimes t^{\alpha_{k+1}})\\
&\quad =\left((a_{k}(z\q^{l})\otimes t^{\alpha_k +r})_{-1} (a_{k+1}(z)\otimes t^{\alpha_{k+1}})\right)\bigg.\bigg|_{z=z_{k+1}\q^{\alpha_1+...+\alpha_{k-1}}}.
\end{align*}
Since polynomials in \eqref{dong_polinomial} do not have any zeroes in the set $\left\{\q^{n}:n\in\mathbb{Z}_{\geq 0} \right\}$, we conclude  that \eqref{dong8} is quasi-commutative. 
\end{prf}

Denote by $R_{z_j}$ the operator $$R_{z_j}\,\colon\,a(z_1,...,z_{j-1},z_j,z_{j+1},...,z_n)\,\mapsto\, a(z_1,...,z_{j-1},z_j \q,z_{j+1},...,z_n).$$
We will need the following three technical lemmas. 

\begin{lem}\label{lema_08} 
For every nonnegative integer $s$ and a polynomial  $p(z)\in \mathbb{F}(\q)[z]_0$ we have
\begin{align}
&\left(\frac{\partial_{\q}}{\partial_{\q} z_1} R_z + \frac{\partial_{\q}}{\partial_{\q} z}\right)^s = 
\sum_{l=0}^{s} \gauss{s}{l}_\q R_{z}^{l}\frac{\partial_{\q}^{l}}{\partial_{\q} z_1^l} 
\frac{\partial_{\q}^{s-l}}{\partial_{\q} z^{s-l}};\label{label1} \\
&R_z \frac{\partial_{\q}}{\partial_{\q} z_1} p(z/z_1)=-\frac{z}{z_1}\frac{\partial_{\q}}{\partial_{\q} z} p(z/z_1).\label{lemma:1111}
\end{align}
\end{lem}

\begin{prf}
Equality \eqref{label1} can be proved by induction over $s$, using
$$\displaystyle\frac{\partial_{\q}}{\partial_{\q} z} R_z = \q R_z \frac{\partial_{\q}}{\partial_{\q} z},\qquad \displaystyle\frac{\partial_{\q}}{\partial_{\q} z_1} R_z=R_z\displaystyle\frac{\partial_{\q}}{\partial_{\q} z_1} $$
and Proposition \ref{expand}; Equality \eqref{lemma:1111} can be verified by a direct calculation.
\end{prf}

\begin{lem}\label{lema_10}
Let $a(z),b(z)\in\mathcal{E}(L)$ and $p(z)\in\mathbb{F}(\q)[z]_0$ such that
\begin{align}\label{assumx}
p(z/z_1)a(z_1)b(z),\,p(z/z_1)a^{[1]}(z_1)b(z\q),\,p(z/z_1)a(z_1)b^{[1]}(z)\in\om(L,L((z_1, z))).
\end{align}
Then
\begin{align}
&\frac{d_{\q}}{d_{\q} z}\lim_{z_1\to z} p(z/z_1)a(z_1)b(z)
=\lim_{z_1\to z} \left(\frac{\partial_{\q}}{\partial_{\q} z_1} R_z + \frac{\partial_{\q}}{\partial_{\q} z}\right)p(z/z_1)a(z_1)b(z);\label{eq:1111}\\
&\frac{d_{\q} }{d_{\q} z}\lim_{z_1\to z} p(z/z_1)a(z_1)b(z)
=\lim_{z_1\to z} p(z/z_1)\left(\frac{\partial_{\q}}{\partial_{\q} z_1} R_z + \frac{\partial_{\q}}{\partial_{\q} z}\right)a(z_1)b(z).\label{eq:1112}
\end{align}
\end{lem}

\begin{prf}
Equality \eqref{eq:1111} is obvious. Let us prove \eqref{eq:1112}. We have
{\allowdisplaybreaks
\begin{align}
&\lim_{z_1\to z} \left(\frac{\partial_{\q}}{\partial_{\q} z_1} R_z + \frac{\partial_{\q}}{\partial_{\q} z}\right)p(z/z_1)a(z_1)b(z)\label{eq:11125}\\
=&\lim_{z_1\to z} R_z \left(
\left(\frac{\partial_{\q}}{\partial_{\q} z_1}p(z/z_1)\right) a(z_1 \q)b(z)
+p(z/z_1)a^{[1]}(z_1)b(z)
\right)\label{eq:1113}\\
&+\lim_{z_1\to z}  \left(
\left(\frac{\partial_{\q}}{\partial_{\q} z}p(z/z_1)\right)a(z_1 )b(z\q)
+p(z/z_1)a(z_1)b^{[1]}(z)
\right).\label{eq:1114}
\end{align}}
Note that \eqref{assumx}  implies that limits \eqref{eq:1113} and \eqref{eq:1114} do exist.
Using  \eqref{lemma:1111} we can express \eqref{eq:1113} as
{\allowdisplaybreaks
\begin{align}
&\lim_{z_1\to z} R_z \left(
\left(\frac{\partial_{\q}}{\partial_{\q} z_1}p(z/z_1)\right)a(z_1 \q)b(z)
+p(z/z_1)a^{[1]}(z_1)b(z)
\right)\nonumber\\
=&\lim_{z_1\to z}  \left(
\left(-\frac{z}{z_1}\cdot\frac{\partial_{\q}}{\partial_{\q} z}p(z/z_1)\right)a(z_1 \q)b(z\q)
+p(z\q /z_1)a^{[1]}(z_1)b(z\q)
\right)\nonumber\\
=&\lim_{z_1\to z}  \left(
\left(-\frac{\partial_{\q}}{\partial_{\q} z}p(z/z_1)\right)a(z_1 \q)b(z\q)
+\frac{z_1}{z}p(z\q /z_1)a^{[1]}(z_1)b(z\q)
\right).\label{eq:1116}
\end{align}}
Combining \eqref{eq:11125} and \eqref{eq:1116} we get
{\allowdisplaybreaks\begin{align}
&\lim_{z_1\to z} \left(\frac{\partial_{\q}}{\partial_{\q} z_1} R_z + \frac{\partial_{\q}}{\partial_{\q} z}\right)p(z/z_1)a(z_1)b(z)\nonumber\\
=&\lim_{z_1\to z}  \left(
\left(-\frac{\partial_{\q}}{\partial_{\q} z}p(z/z_1)\right)a(z_1 \q)b(z\q)
+\frac{z_1}{z}p(z\q /z_1)a^{[1]}(z_1)b(z\q)
\right.\nonumber\\
&\qquad + \left.
\left(\frac{\partial_{\q}}{\partial_{\q} z}p(z/z_1)\right)a(z_1 )b(z\q)
+p(z/z_1)a(z_1)b^{[1]}(z)
\right)\nonumber\\
=&\lim_{z_1\to z}\left( \left(
-\frac{\partial_{\q}}{\partial_{\q} z}p(z/z_1)\right)a^{[1]}(z_1)b(z\q)z(\q-1)\right.\nonumber\\
&\qquad+\Bigg.\frac{z_1}{z}p(z\q /z_1)a^{[1]}(z_1)b(z\q)+p(z/z_1)a(z_1)b^{[1]}(z)
\Bigg)\nonumber\\
=&\lim_{z_1\to z}\left( \left(\left(
-\frac{\partial_{\q}}{\partial_{\q} z}p(z/z_1)\right)z(\q-1)+\frac{z_1}{z}p(z\q /z_1)\right)a^{[1]}(z_1)b(z\q)+p(z/z_1)a(z_1)b^{[1]}(z)
\right)\nonumber\\
=&\lim_{z_1\to z}\left( \left(\left(
-\frac{\partial_{\q}}{\partial_{\q} z}p(z/z_1)\right)z(\q-1)+\frac{z_1}{z}p(z\q /z_1)\right)p(z/z_1)a^{[1]}(z_1)b(z\q)\right.\nonumber\\
&\qquad+\Bigg. p(z/z_1)^2 a(z_1)b^{[1]}(z)
\Bigg).\label{eq:1119}
\end{align}}
Since
$$\lim_{z_1\to z} \left(\left(
-\frac{\partial_{\q}}{\partial_{\q} z}p(z/z_1)\right)z(\q-1)+\frac{z_1}{z}p(z\q /z_1)\right)=1,$$
Assumption \eqref{assumx}, together with \eqref{eq:1119}, implies
\begin{align*}
&\lim_{z_1\to z} \left(\frac{\partial_{\q}}{\partial_{\q} z_1} R_z + \frac{\partial_{\q}}{\partial_{\q} z}\right)p(z/z_1)a(z_1)b(z)\\
&\quad =\lim_{z_1\to z}\left( p(z/z_1)a^{[1]}(z_1)b(z\q)+p(z/z_1) a(z_1)b^{[1]}(z)
\right).
\end{align*}
Finally, Equality \eqref{eq:1112} follows from \eqref{eq:1111}.
\end{prf}

Using  \eqref{label1}  and Lemma \ref{lema_10} one can  prove:

\begin{lem}\label{anotherone}
Let $(a (z)\otimes t^{\alpha},b (z)\otimes t^{\beta})$ be a quasi-commutative pair in $\mathcal{E}(L)_t$,  $r,s$ nonnegative integers and $p(z)\in\mathbb{F}(\q)[z]_0$ a polynomial satisfying 
\begin{align*}
&p(z/z_1)a(z)b(z\q^{\alpha})\in \om(L,L((z_1,z)));\\
&p(z/z_1)a^{[r+l]}(z_1)b^{[s-l]}(z\q^{\alpha+l +r})\in \om(L,L((z_1,z)))\quad \text{for }l=0,1,...,s.
\end{align*}
Then
\begin{align*}
&\frac{\partial_\q^s}{\partial_{\q} z_{1}^{s}}\left(\lim_{z_3 \to z_1}p(z_1/z_3) a^{[r]}(z_3) b(z_1\q^{r+\alpha})\right)\\
&\quad  =\sum_{l=0}^{s} 
\gauss{s}{l}_{\q}\q^{(s-l)(r+\alpha )}\lim_{z_3 \to z_1} p(z_1 / z_3) a^{[r+l]}(z_3)b^{[s-l]}(z_1 \q^{\alpha+l +r}).
\end{align*}
\end{lem}

The set $\mathcal{S}\subseteq\mathcal{E}(L)_t$ is said to be {\em quasi-commutative} if every finite sequence in $\mathcal{S}$ is quasi-commutative.
The quasi-commutative set $\mathcal{S}\subseteq\mathcal{E}(L)_t$ is said to be {\em closed} if $a(z,t)_r b(z,t)$ is an element of $\mathcal{S}$ for all
$a(z,t),b(z,t)\in\mathcal{S}$, $r\in\mathbb{Z}$.

\begin{lem}\label{main}
Let $V$ be a closed, quasi-commutative subspace of $\mathcal{E}(L)_t$ such that $\vac\in V$. Then $(V,Y,\vac)$ satisfies all the axioms
of the graded nonlocal $\q$-vertex algebra, except maybe grading restrictions \eqref{VVV} and \eqref{V4}. 
\end{lem}

\begin{prf}
The definition of the operator $Y$ is given in \eqref{sydney}, so \eqref{V0} is satisfied. Definition \ref{produkti} implies \eqref{V1}.
The definition of the vector $\vac\in\mathcal{E}(L)_t$ is given in \eqref{vacuum} and  Corollary \ref{lotsofthings} implies \eqref{V2} and \eqref{V3}.

Let us prove \eqref{V5}.
It is sufficient to consider an arbitrary quasi-commutative sequence of homogeneous elements $(a(z)\otimes t^{\alpha},b(z)\otimes t^{\beta},c(z)\otimes t^{\gamma})$, whose all subsequences are also quasi-commutative.
Fix nonnegative integers $r_0,s_0$. Let 
$p(z)\in\mathbb{F}(\q)[z]_0$ be a polynomial satisfying 
{\allowdisplaybreaks\begin{align}
&p(z/z_1) a^{[r]}(z_1) b(z\q^{r+\alpha})\in\om(L,L((z_1,z)));\label{as1}\\
&p(z_1 / z_3) a^{[r+l]}(z_3)b^{[s-l]}(z_1 \q^{\alpha+l+r})\in\om(L,L((z_3,z_1)));\label{as2}\\
& p(z/z_1)p(z/z_3)p(z_1 / z_3) a^{[r+l]}(z_3)b^{[s-l]}(z_1 \q^{\alpha+l+r})  c(z\q^{\alpha+\beta+r+s+k})\in\om(L,L((z_3,z_1,z)));\label{as3}\\
& p(z/z_1) b^{[v]}(z_1)c(z\q^{\beta +v})\in\om(L,L((z_1,z)))\label{as4}
\end{align}}
for all 
$$r=0,1,...,r_0+s_0, \quad k,s,v=0,1,...,s_0,\quad l=0,1...,s. $$  
Using  
\eqref{as3},  \eqref{as4} and Lemma \ref{mp}
we get
{\allowdisplaybreaks
\begin{align}
&\lim_{z_1,z_3 \to z} p(z/z_1)p(z/z_3 ) p(z_1 /z_3) a^{[u]}(z_3)
 b^{[v]}(z_1\q^{\alpha+u})c(z\q^{\alpha+\beta +u+v})\nonumber\\
 =&\lim_{z_3 \to z}\lim_{z_1 \to z} p(z/z_1)p(z/z_3 ) p(z_1 /z_3) a^{[u]}(z_3)
  b^{[v]}(z_1\q^{\alpha+u})c(z\q^{\alpha+\beta +u+v})\nonumber\\
   =&\lim_{z_3 \to z}p(z/z_3 )^2 a^{[u]}(z_3)\lim_{z_1 \to z} p(z/z_1)  
    b^{[v]}(z_1\q^{\alpha+u})c(z\q^{\alpha+\beta +u+v})\label{as5}
\end{align}}
for all $u=0,1,...,r_0 +s_0$, $v=0,1,...,s_0$.
Using   \eqref{as2}, \eqref{as3} and Lemma \ref{mp} we get
{\allowdisplaybreaks
\begin{align}
&\lim_{z_1,z_3 \to z} p(z/z_1)p(z/z_3)p(z_1 / z_3) a^{[r+l]}(z_3)b^{[s-l]}(z_1 \q^{\alpha+l +r})  c(z\q^{\alpha+\beta+r+s})\nonumber\\
=& \lim_{z_1 \to z}\lim_{z,z_3 \to z_1} p(z/z_1)p(z/z_3)p(z_1 / z_3) a^{[r+l]}(z_3)b^{[s-l]}(z_1 \q^{\alpha+l +r})  c(z\q^{\alpha+\beta+r+s})\nonumber\\
=& \lim_{z_1 \to z}\lim_{z \to z_1}\lim_{z_3 \to z_1} p(z/z_1)p(z/z_3)p(z_1 / z_3) a^{[r+l]}(z_3)b^{[s-l]}(z_1 \q^{\alpha+l +r})  c(z\q^{\alpha+\beta+r+s})\nonumber\\
=& \lim_{z_1 \to z}\lim_{z \to z_1} p(z/z_1)^2\lim_{z_3 \to z_1}\left(p(z_1 / z_3) a^{[r+l]}(z_3)b^{[s-l]}(z_1 \q^{\alpha+l +r})\right)  c(z\q^{\alpha+\beta+r+s})\nonumber\\
=& \lim_{z_1 \to z} p(z/z_1)^2\lim_{z_3 \to z_1}\left(p(z_1 / z_3) a^{[r+l]}(z_3)b^{[s-l]}(z_1 \q^{\alpha+l +r})\right)  c(z\q^{\alpha+\beta+r+s})\label{as6}
\end{align}}
for all $r=0,1,...,r_0$, $s=0,1,...,s_0$, $l=0,1,...,s$.

For  Taylor  series $A(z_0,z_2), B(z_0,z_2)\in\mathcal{E}(L) [[z_0,z_2]]_t$ we write $A(z_0,z_2)\equiv B(z_0,z_2)$ if  the coefficient of $z_0^r z_2^s$ in $A(z_0,z_2)- B(z_0,z_2)$ equals zero for all 
$r\leq r_0$ and $s\leq s_0$. 
In the following calculations we assume that the variables $z_0,z_2$ satisfy \eqref{commcons}.
We use the relation "$\equiv$"
because the polynomial $p(z)$ satisfying \eqref{as1}--\eqref{as4} depends on the choice of nonnegative integers $r_0,s_0$.
Using \eqref{as1}, \eqref{as6} and Lemma \ref{anotherone} we get
{\allowdisplaybreaks
\begin{align}
&Y(Y(a(z)\otimes t^{\alpha} ,z_0)b(z)\otimes t^{\beta},z_2)(c(z)\otimes t^{\gamma})\nonumber\\
\equiv &Y( \sum_{r=0}^{r_0}\lim_{z_1 \to z}\frac{p(z/z_1)}{[r]_{\q}!} a^{[r]}(z_1) b(z\q^{r+\alpha})\otimes t^{\alpha+\beta+r}z_{0}^{r} ,z_2 )(c(z)\otimes t^{\gamma})\nonumber\\
=& \sum_{r=0}^{r_0}z_{0}^{r}Y(\lim_{z_1 \to z}\frac{p(z/z_1)}{[r]_{\q}!} a^{[r]}(z_1) b(z\q^{r+\alpha})\otimes t^{\alpha+\beta+r} ,z_2 )(c(z)\otimes t^{\gamma})\nonumber\\
\equiv& \sum_{r=0}^{r_0}\sum_{s=0}^{s_0}\lim_{z_1 \to z}\left(\frac{p(z/z_1)^2}{[s]_{\q}!}\left(\lim_{z_3 \to z_1}\frac{p(z_1/z_3)}{[r]_{\q}!} a^{[r]}(z_3) b(z_1\q^{r+\alpha})\right)^{[s]}  c(z\q^{\alpha+\beta+r+s})\right)\nonumber\\
&\otimes t^{\alpha+\beta+\gamma+r+s}z_{0}^{r}z_{2}^s\nonumber\\
=& \sum_{r=0}^{r_0}\sum_{s=0}^{s_0}\sum_{l=0}^{s} 
\frac{\q^{(s-l)(\alpha +r)}}{[r]_\q ! [l]_\q ! [s-l]_\q !}\nonumber\\
&\cdot\lim_{z_1 \to z}\left( p(z/z_1)^2\left(\lim_{z_3 \to z_1}p(z_1 / z_3) a^{[r+l]}(z_3)b^{[s-l]}(z_1 \q^{\alpha+l +r})\right)  c(z\q^{\alpha+\beta+r+s})\right)\nonumber\\
&\quad\otimes t^{\alpha+\beta+\gamma+r+s}z_{0}^{r}z_{2}^s\nonumber\\
=& \sum_{r=0}^{r_0}\sum_{s=0}^{s_0}\sum_{l=0}^{s} 
\frac{\q^{(s-l)(\alpha +r)}}{[r]_\q ! [l]_\q ! [s-l]_\q !}\nonumber\\
&\cdot\hspace{-6pt}\lim_{z_1,z_3 \to z}\left( p(z/z_1)p(z/z_3)p(z_1 / z_3) a^{[r+l]}(z_3)b^{[s-l]}(z_1 \q^{\alpha+l +r})  c(z\q^{\alpha+\beta+r+s})\right)\otimes t^{\alpha+\beta+\gamma+r+s}z_{0}^{r}z_{2}^s.\label{number1}
\end{align}}
Using \eqref{as4} and \eqref{as5} we get
{\allowdisplaybreaks
\begin{align}
&Y(a(z)\otimes t^{\alpha},z_0 +z_2)Y(b(z)\otimes t^{\beta},z_2 \q^{\wt a(z)\otimes t^{\alpha}}) (c(z)\otimes t^{\gamma})\nonumber\\
\equiv&Y(a(z)\otimes t^{\alpha},z_0 +z_2)\sum_{v=0}^{s_0}
\lim_{z_1 \to z} \frac{p(z/z_1)\q^{\alpha v}}{[v]_\q !} b^{[v]}(z_1)c(z\q^{\beta +v})\otimes t^{\beta +\gamma+v} z_2^{v}\nonumber\\
\equiv&\sum_{u=0}^{r_0 +s_0}\sum_{v=0}^{s_0}
\lim_{z_3 \to z}\left(\frac{p(z/z_3 )^2}{[u]_\q !} a^{[u]}(z_3)
\left(\lim_{z_1 \to z} \frac{p(z/z_1)\q^{\alpha v}}{[v]_\q !} b^{[v]}(z_1\q^{u+\alpha})c(z\q^{\alpha+\beta +u+v})\right)\right)\nonumber\\
&\otimes t^{\alpha+\beta+\gamma +u+v} z_2^{v}(z_0 +z_2)^u\nonumber\\
=&\sum_{u=0}^{r_0 +s_0}\sum_{v=0}^{s_0}
\lim_{z_3 \to z}\left(\frac{p(z/z_3 )^2}{[u]_\q !} a^{[u]}(z_3)
\left(\lim_{z_1 \to z} \frac{p(z/z_1)\q^{\alpha v}}{[v]_\q !} b^{[v]}(z_1\q^{u+\alpha})c(z\q^{\alpha+\beta +u+v})\right)\right)\nonumber\\
&\otimes t^{\alpha+\beta+\gamma +u+v} z_2^{v}\left(\sum_{l=0}^{u}\gauss{u}{l}_{\q}z_0^{u-l}
z_2^l\right)\nonumber\\
=&\sum_{u=0}^{r_0 +s_0}\sum_{v=0}^{s_0}\sum_{l=0}^{u}
\frac{\q^{v(\alpha+u-l) }}{[u-l]_\q ! [l]_\q ! [v]_\q !}\nonumber\\
&\cdot
\lim_{z_3 \to z}\left(p(z/z_3 )^2 a^{[u]}(z_3)
\left(\lim_{z_1 \to z} p(z/z_1) b^{[v]}(z_1\q^{u+\alpha})c(z\q^{\alpha+\beta +u+v})\right)\right)\otimes t^{\alpha+\beta +\gamma+u+v} z_0^{u-l}z_2^{t+l}\nonumber\\
=&\sum_{u=0}^{r_0 +s_0}\sum_{v=0}^{s_0}\sum_{l=0}^{u}
\frac{\q^{v(\alpha+u-l) }}{[u-l]_\q ! [l]_\q ! [v]_\q !}\nonumber\\
&\cdot
\lim_{z_1,z_3 \to z}\left(p(z/z_1)p(z/z_3 ) p(z_1 /z_3) a^{[u]}(z_3)
 b^{[v]}(z_1\q^{u+\alpha})c(z\q^{\alpha+\beta +u+v})\right)\otimes t^{\alpha+\beta +\gamma+u+v} z_0^{u-l}z_2^{v+l}.\label{number2}
\end{align}}

Applying the substitutions $u=r+l$ and $v=s-l$ to \eqref{number2} we get $\text{\eqref{number1}}\equiv\text{\eqref{number2}}$, i.e.
$$Y(Y(a(z)\otimes t^{\alpha} ,z_0)b(z)\otimes t^{\beta},z_2)c(z)\otimes t^{\gamma}
\equiv
Y(a(z)\otimes t^{\alpha},z_0 +z_2)Y(b(z)\otimes t^{\beta},z_2 \q^{\alpha}) c(z)\otimes t^{\gamma}.
$$
Finally, since $r_0$ and $s_0$ were arbitrary nonnegative integers, we conclude
$$Y(Y(a(z)\otimes t^{\alpha} ,z_0)b(z)\otimes t^{\beta},z_2)c(z)\otimes t^{\gamma}
=
Y(a(z)\otimes t^{\alpha},z_0 +z_2)Y(b(z)\otimes t^{\beta},z_2 \q^{\alpha}) c(z)\otimes t^{\gamma}.
$$
\end{prf}

\begin{lem}\label{important}
Let $V$ be a maximal quasi-commutative subspace of $\mathcal{E}(L)_t$. Then
\begin{enumerate}[(a)]
\item\label{473829471} $V$ is closed.
\item\label{473829472} $V$ contains $\vac$.
\item\label{473829473} If $\sum_{i=1}^{N}a_{i}(z)\otimes t^{\alpha_i} \in V$, then $a_{i}(z)\otimes t^{\alpha_i}\in V$ for $i=1,2,...,N$.
\item\label{473829474} $V$ is a graded nonlocal $\q$-vertex algebra.
\end{enumerate}
\end{lem}

\begin{prf}
Statement \eqref{473829471} follows from Lemma \ref{dong}. Maximality of $V$, together with Corollary \ref{lotsofthings}, implies 
\eqref{473829472}. Statement \eqref{473829473} is a consequence of the maximality of $V$ and Definition \ref{qkomm}.
Finally, let $V_{(u)}$ be a subspace of $V$ consisting of all homogeneous elements $a(z)\otimes t^{u} \in V$.
Lemma \ref{main}, Statements \eqref{473829471}--\eqref{473829473} and Corollary \ref{lotsofthings} imply \eqref{473829474}.
\end{prf}

The following theorem is the main result in this section.

\begin{thm}\label{main2}
Let $\mathcal{S}$ be a quasi-commutative subset of $\mathcal{E}(L)_t$. There exists a unique
smallest graded nonlocal q-vertex algebra $V\subseteq\mathcal{E}(L)_t$ such that $\mathcal{S}\subseteq V$.
Furthermore, if $\mathcal{S}$ consists of homogeneous elements, we have
\begin{equation}\label{span}
V=\xpan \left\{a_1(z,t)_{r_1}\cdots a_k(z,t)_{r_k}\vac\,:\,a_j(z,t)\in\mathcal{S},\,r_j< 0,\,j=1,...,k,\,k\in\mathbb{Z}_{\geq 0} \right\}.
\end{equation}
\end{thm}

\begin{prf}
By Zorn's Lemma, $\mathcal{S}$ is a subset of the maximal closed quasi-commutative subset of $\mathcal{E}(L)_t$ which is, by Lemma \ref{important},
graded nonlocal $\q$-vertex algebra. Therefore, we can construct $V$ as the intersection of all graded nonlocal $\q$-vertex algebras  containing $\mathcal{S}$, whose grading is defined as in the proof of Lemma \ref{important}.

Let us prove \eqref{span}. Denote the right-hand side in \eqref{span} by $R$. It is evident that $R\subseteq V$. By \eqref{V5} (cf. also \eqref{asocc})
$R$ is closed and by Lemma \ref{dong} $R$ is quasi-commutative. Therefore, since $\vac\in R$, by Lemma \ref{main}, $R$ is  graded nonlocal $\q$-vertex algebra
(notice that grading restrictions \eqref{VVV} and \eqref{V4} hold because $\mathcal{S}$ consists of homogeneous elements). 
By construction, $V$ is the smallest graded nonlocal $\q$-vertex algebra containing $\mathcal{S}$, so $V\subseteq R$.
\end{prf}

\begin{rem}\label{algebra}
By taking $r=s=-1$ in \eqref{asocc} we get
$(a_{-1}b)_{-1}c = a_{-1}(b_{-1} c)$
for vertex operator products given by Definition \ref{produkti}.
Hence the graded nonlocal $\q$-vertex algebra $V$ constructed in  Theorem \ref{main2} becomes an unital associative algebra  with product being defined by
$a\cdot b=a_{-1} b$ for all $a,b\in V$.
\end{rem}

\subsection{Examples}\label{examplesection}

Consider the quantum affine algebra $U_q(\widehat{\mathfrak{g}})$ given by Definition \ref{drinfeld}, when $\widehat{\mathfrak{g}}$ is an affine Kac-Moody
Lie algebra of type $(ADE)^{(1)}$ and $\widehat{A}=(a_{ij})_{i,j=0}^{n}$ is its Cartan matrix. Let $L$ be a restricted $U_q(\widehat{\mathfrak{g}})$-module. Set
$$\mathcal{S}=\left\{x_{i}^{+}(z)\otimes 1\,:\,i=1,2,...,n\right\}\subset \mathcal{E}(L)_t\quad\text{and}\quad \q=q^2.$$
Recall definitions of the vertex operator map $Y$ and the vacuum vector $\vac$ in \eqref{sydney} and \eqref{vacuum}.

\begin{thm}\label{construction}
There exists a unique smallest graded nonlocal $\q$-vertex algebra $(\left<\mathcal{S}\right>,Y,\vac)$ which contains $\mathcal{S}$.
Furthermore, $\left<\mathcal{S}\right>$ is spanned by
\begin{equation*}
 \left\{(x_{i_1}^{+}(z)\otimes 1)_{r_1}\cdots (x_{i_m}^{+}(z)\otimes 1)_{r_m}\vac\,:\,1\leq i_j\leq n,\,r_j< 0,\,j=1,...,m,\,m\in\mathbb{Z}_{\geq 0} \right\}.
\end{equation*}
\end{thm}

\begin{prf}
It is sufficient to prove that every finite sequence in $\mathcal{S}$ is quasi-commutative. Then,  we can apply Theorem \ref{main2}. 
First, by considering relations \eqref{D8} and \eqref{Da} we see that every pair in $\mathcal{S}$ is quasi-commutative. 
Note that for $(\alpha_i,\alpha_j)=0$ relation \eqref{D8}  gives us
$$ (z_1-z_2)x_{i}^{+}(z_1)x_{j}^{+}(z_2)=(z_1-z_2)x_{j}^{+}(z_2)x_{i}^{+}(z_1),$$
so we have to use \eqref{Da} to establish quasi-commutativity.

The rest of the proof goes analogously as the proof of \cite[Lemma 3.2]{Li1}. Assume that every sequence of length $m$ in $\mathcal{S}$ is
quasi-commutative.  We shall prove that every sequence of length $m+1$ in $\mathcal{S}$ is quasi-commutative, so the statement will follow by induction.  Consider the sequence $(a_1(z),...,a_{m+1}(z))$ in $\mathcal{S}$. Let $p_1(z),p_2(z),p_3(z)\in\mathbb{C}(q^{1/2})[z]_0$, $p_4(z)\in\mathbb{C}(q^{1/2})[z]$ be polynomials satisfying
\begin{align}
&p_1(z_2 /z_1) a_1 (z_1) a_2 (z_2)=p_4(z_2/z_1)  a_2 (z_2)a_1 (z_1) \in \om(L,L((z_1,z_2)));\label{ll1}\\
&P_2(z_2,...,z_{m+1}) a_2 (z_2)\cdots a_{m+1} (z_{m+1}) \in \om(L,L((z_2,...,z_{m+1})));\label{ll2}\\
&P_3(z_1,z_3,...,z_{m+1})a_1 (z_1)a_3 (z_3)\cdots a_{m+1} (z_{m+1}) \in \om(L,L((z_1,z_3,...,z_{m+1}))),\label{ll3}
\end{align}
where
\begin{align*}
P_2(z_2,...,z_{m+1})=\prod_{2\leq i< j\leq m+1}p_2 (z_j /z_i),\quad P_3(z_1,z_3,...,z_{m+1})=\prod_{\substack{1\leq i< j\leq m+1\\i,j\neq 3}}p_3 (z_j /z_i).
\end{align*}
Of course, existence of polynomials $p_1(z)$ and $ p_4(z)$ satisfying \eqref{ll1} is a consequence of \eqref{D8} and \eqref{Da}.
Using \eqref{ll1}--\eqref{ll3} we get
\begin{align}
&p_1(z_2 /z_1)P_2(z_2,...,z_{m+1})P_3(z_1,z_3,...,z_{m+1})a_1 (z_1)a_2 (z_2)a_3 (z_3)\cdots a_{m+1} (z_{m+1})\label{ll4}\\
&=p_4(z_2 /z_1)P_2(z_2,...,z_{m+1})P_3(z_1,z_3,...,z_{m+1})a_2 (z_2)a_1 (z_1)a_3 (z_3)\cdots a_{m+1} (z_{m+1}).\label{ll5}
\end{align}
Since \eqref{ll4} is an element of $\om(L,L((z_1))((z_2,z_3,...,z_{m+1})))$ and  \eqref{ll5} is an element of $\om(L,L((z_2))((z_1,z_3,...,z_{m+1})))$,
we conclude that the both sides are contained in $\om(L,L((z_1,z_2,z_3,...,z_{m+1})))$. Since $p_1(z), p_2(z), p_3(z)\in\mathbb{C}(q^{1/2})[z]_0$,  the statement follows by induction.
\end{prf}

Naturally, there is an analogous result for negative current operators. Set 
$$\mathcal{T}=\left\{x_{i}^{-}(z)\otimes 1\,:\,i=1,2,...,n\right\}\subset \mathcal{E}(L)_t\quad\text{and}\quad \q=q^{-2}.$$

\begin{kor}
There exists a unique smallest graded nonlocal $\q$-vertex algebra $(\left<\mathcal{T}\right>,Y,\vac)$ which contains $\mathcal{T}$.
Furthermore, $\left<\mathcal{T}\right>$ is spanned by
\begin{equation*}
\left\{(x_{i_1}^{-}(z)\otimes 1)_{r_1}\cdots (x_{i_m}^{-}(z)\otimes 1)_{r_m}\vac\,:\,1\leq i_j\leq n,\,r_j< 0,\,j=1,...,m,\,m\in\mathbb{Z}_{\geq 0} \right\}.
\end{equation*}
\end{kor}

By repeating the proof of Theorem \ref{construction} almost verbatim, one can prove:

\begin{kor}\label{cor}
For every $(\alpha_1,...,\alpha_n)$ in $\left(\mathbb{Z}_{\geq 0}\right)^{\times n}$ there exists a unique smallest graded nonlocal $\q$-vertex algebra $(V,Y,\vac)$, $\q=q^{\pm 2}$, generated by the set
$$\left\{x_{i}^{\pm}(z)\otimes t^{\alpha_i}\,:\,i=1,2,...,n\right\}\subset \mathcal{E}(L)_t.$$
Furthermore, $V$ is spanned by
\begin{equation*}
\left\{(x_{i_1}^{\pm}(z)\otimes t^{\alpha_{i_1}})_{r_1}\cdots (x_{i_m}^{\pm}(z)\otimes t^{\alpha_{i_m}})_{r_m}\vac\,:\,1\leq i_j\leq n,\,r_j< 0,\,j=1,...,m,\,m\in\mathbb{Z}_{\geq 0} \right\}.
\end{equation*}
\end{kor}

\begin{rem}
In this paper, we are mainly concerned with the vertex algebraic structures arising from $U_q(\widehat{\mathfrak{sl}}_2)$. 
Even though  the construction of a much broader class of graded nonlocal $\q$-vertex algebras was here presented, 
the role of these examples was to clarify the application of Theorem \ref{main2}. 
The case $\widehat{\mathfrak{g}}=\widehat{\mathfrak{sl}}_2$, which is studied in the next section,
gives rise to the structures whose character formulae  greatly resemble those from the classical theory.
However, the definition of vertex operator products for  $\widehat{\mathfrak{g}}\neq \widehat{\mathfrak{sl}}_2$,
which would provide the analogous correspondence, might require the use of bigger spaces $\mathcal{E}(L)\otimes\mathbb{F}[[t_1,...,t_n]]$ in a similar vein as in
\cite{Kozic}. 
\end{rem}

\section{Graded nonlocal \texorpdfstring{$\q$}{q}-vertex algebras for \texorpdfstring{$U_q(\widehat{\mathfrak{sl}}_2)$}{Uq(sl2)}}\label{section3}

\subsection{Quasi-particles}

In this section, we consider quantum affine algebra $U_q(\widehat{\mathfrak{sl}}_2)$.
Let $c\geq 1$ and
$$L=(L_{0})^{\otimes c}=(M(1)\otimes \mathbb{C}\left\{Q\right\})^{\otimes c}.$$
In the rest of this paper, we will denote the operator  $x_{1}^{+}(z)\in\mathcal{E}(L)$ by $x(z)$ and the operator
$\phi_1 (z)\in\mathcal{E}(L)$ by $\phi(z)$.
By Corollary \ref{cor}  for every $n\geq 0$ the operator $x(z)\otimes t^{n}\in\mathcal{E}(L)_t$
generates a  graded nonlocal $\q$-vertex algebra, which will be denoted  by $V_{c,n}$.
Our goal is to construct  monomial bases for $V_{c,1}$ and, consequently, obtain corresponding
character formulae. The main building blocks of basis elements will be quasi-particles.

\begin{defn}\label{definition}
For a positive integer $m$ we call the operator
$$x_{m,n} (z,t)=\underbrace{(x(z)\otimes t^{n})_{-1}(x(z)\otimes t^{n})_{-1}...(x(z)\otimes t^{n})_{-1}}_{m} \vac \in V_{c,n}$$
a {\em quasi-particle of charge} $m$.
\end{defn}

\begin{pro}\label{mkmk}
For all $m,k\geq 1$, $n\geq 0$ we have in $V_{c,n}$
\begin{equation}\label{mkmkmk}
x_{m,n}(z,t)_{-1}x_{k,n}(z,t) = x_{m+k,n}(z,t).
\end{equation}
\end{pro}

\begin{prf}
The statement  is  a consequence of Remark \ref{algebra}.
\end{prf}

\begin{rem}
Definition \ref{definition} and Equality \eqref{mkmkmk} coincide with the
properties of the original quasi-particle operators for the affine Kac-Moody Lie algebras (cf. \cite{FS},\cite{G},\cite{Bu2}).
Also, in \cite{Ko1}, (quantum) quasi-particle operators for quantum affine algebra $U_q (\widehat{\mathfrak{sl}}_2)$,
\begin{equation}\label{quasi-quantum}
x^{+}_{m\alpha_{1}}(z):=\lim_{z_{j}\to zq^{2(j-1)}}\left(\prod_{r=1}^{m-1}\prod_{s=r+1}^{m}
\left(1-q^{2}\frac{z_{s}}{z_{r}}\right)\right)x(z_{1})\ldots x(z_{m})\in\mathcal{E}(L),
\end{equation}
were considered. Using quantum integrability \eqref{ding_miwa} it was proved that 
\begin{equation}\label{qint5}
x^{+}_{m\alpha_{1}}(z)=0 \quad\text{on any level }c<m\text{ integrable highest weight module.}
\end{equation}
\end{rem}

From now on, we shall consider only quasi-particles $x_{m,1}(z,t)\in V_{c,1}$. In order to simplify our notation we  will denote 
$x_{m,1}(z,t)$ by $x_{m}(z,t)$.

\begin{pro}\label{integrability}
For   $m\geq 1$ we have
$$x_{m} (z,t)=0\text{ in }V_{c,1}\quad\text{if and only if}\quad m>c.$$
\end{pro}

\begin{prf}
Denote by $x_m (z)\in \mathcal{E}(L)$ the left tensor factor of the element in $\mathcal{E}(L)_t=\mathcal{E}(L)\otimes\mathbb{C}(q^{1/2})[t]$ which is obtained by  setting $t=1$ in $x_m (z,t)$. More precisely,
$$x_m (z) \otimes 1 = \left(x_m(z,t)\right)\big|_{t=1}.$$
The operator $x_m(z)$ is equal, up to a nonzero multiplicative factor, to the operator $x_{m\alpha_1}^{+}(z)$ in \eqref{quasi-quantum}. Furthermore, 
\eqref{qint5} (see \cite[Proposition 16]{Ko1}) implies that $x_{m} (z,t)=0$ if $m>c$.
If $m\leq c$, \cite[Theorem 38]{Ko1} provides   a basis for a certain subspace of $L$. The basis consists of the monomials
of the operator's $x_m (z)$ coefficients acting on a certain vector $v\in L$. This  implies that the operator $x_m(z)$ has some nonzero coefficients,  hence
$x_{m} (z,t)\neq 0$.
\end{prf}

\subsection{Basis for \texorpdfstring{$V_{c,1}$}{Vc1}}\label{x:2}
 Corollary \ref{cor} implies

\begin{pro}The set 
\begin{equation*}
\mathcal{S}_{c,1}= \left\{x_1(z,t)_{r_1}x_1(z,t)_{r_2}\ldots x_1(z,t)_{r_k}\vac\,:\,r_j<0,\,j=1,...,k,\,k\in\mathbb{Z}_{\geq 0} \right\}.
\end{equation*}
spans $V_{c,1}$.
\end{pro}

By setting $r=-1$ in \eqref{asocc} we get
\begin{equation}\label{asoccV1}
(a_{-1} b)_s c =\sum_{l\geq 0}\q^{-u(s+l+1)}a_{-1-l}(b_{s+l}c)\quad\text{for all }s\in\mathbb{Z},
\end{equation}
where $a\in\left(V_{c,1}\right)_{(u)}$ for some integer $u\geq 0$ and $b,c\in V_{c,1}$. Rewriting \eqref{asoccV1} we get
\begin{equation}\label{asoccV2}
\q^{-u(s+1)}a_{-1}(b_{s}c) =-\sum_{l\geq 1}\q^{-u(s+l+1)}a_{-1-l}(b_{s+l}c) + (a_{-1} b)_s c. 
\end{equation}
Formula \eqref{asoccV2} and Proposition \ref{mkmk} imply that an    element $x_1(z,t)_{r_1}\ldots x_1(z,t)_{r_k}\vac$
of $\mathcal{S}_{c,1}$ satisfying $r_1,...,r_{j-1},r_{j+1}<-1$ and $r_{j}=-1$ for some $j< k$ can be written as a linear combination of
the summands 
$$x_1(z,t)_{r_1}\ldots x_{1}(z,t)_{r_{j-1}} x_1(z,t)_{s_j}\ldots x_1(z,t)_{s_k} \vac,\quad\text{for some } s_j,...s_k\in\mathbb{Z},\,s_{j}<-1,$$
and the summand
$$x_1(z,t)_{r_1}\ldots x_{1}(z,t)_{r_{j-1}}x_{2}(z,t)_{r_{j+1}}x_{1}(z,t)_{r_{j+2}}\ldots x_1(z,t)_{r_k}\vac.$$
Applying such a procedure an appropriate number of times on  every element of $\mathcal{S}_{c,1}$
and using Proposition \ref{mkmk}, we can construct a spanning set
for $V_{c,1}$, which consists of the elements
\begin{align}\label{back}
x_{m_1}(z,t)_{r_1}x_{m_2}(z,t)_{r_2}\ldots x_{m_k}(z,t)_{r_k}\vac\quad\text{such that}\quad
\text{if }r_{j}=-1\text{ then } r_{j+1}=...=r_{k}=-1,
\end{align}
where $r_j< 0$, $m_j\geq 1$, $j=1,...,k$, $k\in\mathbb{Z}_{\geq 0}$.  
Finally,  Propositions \ref{mkmk} and \ref{integrability} imply 

\begin{pro}\label{spanningset7}
The set
\begin{align}
\mathcal{B}_{c,1}=\big\{&x_{m_1}(z,t)_{r_1}x_{m_2}(z,t)_{r_2}\ldots x_{m_k}(z,t)_{r_k}\vac\,:\big.\nonumber\\
&\big. r_1,...,r_{k-1}\leq -2,\,r_k\leq -1,\,1\leq m_j\leq c,\,j=1,...,k,\,k\in\mathbb{Z}_{\geq 0} \big\}\label{basis785}
\end{align}
spans $V_{c,1}$.
\end{pro}

Our next goal is to prove that the set $\mathcal{B}_{c,1}$ is linearly independent.
The space $L_0=M(1)\otimes \mathbb{C}\left\{Q\right\}$ is an $U_q (\widehat{\mathfrak{sl}}_2)$-module, so
the space $L=(L_0)^{\otimes c}$
has also a $U_q (\widehat{\mathfrak{sl}}_2)$-module structure given by Drinfeld's coproduct formula, 
as explained in Section \ref{pre:1:2}.
For an element
\begin{equation}\label{notethattheoperator}
a(z,t)=x_{m_1}(z,t)_{r_1}x_{m_2}(z,t)_{r_2}\ldots x_{m_k}(z,t)_{r_k}\vac \in \mathcal{B}_{c,1}
\end{equation}
define a $c$-tuple
$$r_{a(z,t)}=(r^{(1)},r^{(2)},...,r^{(c)}),$$
where $r^{(j)}$ denotes the number of indices $i=1,2,...,k$ such that $m_i\geq j$.
Define the following subspaces of $L_0$:
$$(L_0)_s =\left\{v\in L_0\,:\, Kv=q^{2s}v \right\},\quad s\in\mathbb{Z}.$$
Then $L_0=\coprod_{s\in\mathbb{Z}} (L_0)_s$ and we have Georgiev's projection (\cite{G}, cf. also \cite{Ko1})
$$\pi_{a(z,t)}\colon (L_0)^{\otimes c} \to (L_0)_{r^{(1)}}\otimes (L_0)_{r^{(2)}}\otimes \cdots \otimes (L_0)_{r^{(c)}}.$$
The mapping $\pi_{a(z,t)}$ can be generalized to the space of the formal Laurent series with  coefficients in
$(L_0)^{\otimes c}$, so we can consider the projection $$\pi_{a(z,t)} \colon (L_0)^{\otimes c}[[z_1^{\pm 1},z_2^{\pm 1},...]] \to \left((L_0)_{r^{(1)}}\otimes \cdots \otimes (L_0)_{r^{(c)}}\right)[[z_1^{\pm 1},z_2^{\pm 1},...]].$$

The operator $a(z,t)$ in \eqref{notethattheoperator} can be expressed as a  linear combination of the  monomials
\begin{align}
&A_{1,l_1,m_1}^{r_1}(z,t)_{-1}A_{2,l_2,m_2}^{r_2}(z,t)_{-1}\cdots A_{k-1,l_{k-1},m_{k-1}}^{r_{k-1}}(z,t)_{-1} A_{k,l_k,m_k}^{r_k}(z,t),\label{azt}
\end{align}
where
\begin{align}
&A_{j,l_j,m_j}^{r_j} (z,t)=(z^{r_j +1}\otimes t^{-r_j -1}) x(z\q^{M_j +l_j},t)_{-1}x(z\q^{M_j +l_j+1},t)_{-1} ... x(z\q^{M_j +l_j+m_j -1},t)_{-1}\vac,\nonumber\\
&l_j=0,1,...,-r_j-1,\,\,  j=1,2,..,k, \,\,  M_{j}=\sum_{i=1}^{j-1}(-r_i -1).\label{Aovi789}
\end{align}
All the summands in a such linear combination have the same
degree $\deg_t$ in variable $t$.
The operator 
\begin{equation}\label{notethattheoperatorb}
A_{1,0,m_1}^{r_1}(z,t)_{-1}A_{2,0,m_2}^{r_2}(z,t)_{-1}\cdots A_{k-1,0,m_{k-1}}^{r_{k-1}}(z,t)_{-1} A_{k,0,m_k}^{r_k}(z,t)
\end{equation}
will be called the {\em leading term} of $a(z,t)$.

\begin{lem}\label{zero5}
The leading term
 of any operator
$a(z,t)\in \mathcal{B}_{c,1}$  is nonzero.
\end{lem}

\begin{prf}
Using the analogous arguments as in the proof of \eqref{as5}
one can show that  leading term \eqref{notethattheoperatorb}
of the operator $a(z,t)$ in \eqref{notethattheoperator}
can be written as
\begin{align}
A(z,t)=&\lim_{z_{j,s}\to z}\bigg(\bigg(\prod_{1\leq i\leq j\leq k}\prod_{\substack{1\leq r\leq m_i\\ 1\leq s\leq m_j\\(i,r)<(j,s)}}p(z_{j,s}/z_{i,r})
\bigg.\bigg) \label{zero2}\\
&\qquad\bigg.\cdot\prod_{j=1}^{k} \left((z^{r_j +1}\otimes t^{-r_j -1}) x(z_{j,1}\q^{M_j },t)x(z_{j,2}\q^{M_j +1},t) ... x(z_{j,m_j}\q^{M_j+m_j -1},t)\right)\bigg)\nonumber
\end{align}
for some polynomial $p(z)\in\mathbb{C}(q^{1/2})[z]_0$.
Applying $A(z,t)$ on the vector $1=1^{\otimes c}\in (L_0)^{\otimes c}=L$ and
 setting $t=1$ we get an element of $L((z))$:
\begin{align}
&\lim_{z_{j,s}\to z}\bigg(\bigg(\prod_{1\leq i\leq j\leq k}\prod_{\substack{1\leq r\leq m_i\\ 1\leq s\leq m_j\\(i,r)<(j,s)}}p(z_{j,s}/z_{i,r})
\,\bigg.\bigg)\nonumber\\
&\qquad\bigg.\cdot\,\prod_{j=1}^{k} \left(z^{r_j +1} x(z_{j,1}\q^{M_j })x(z_{j,2}\q^{M_j +1}) ... x(z_{j,m_j}\q^{M_j+m_j -1})\right)\bigg)\,1.\label{complicated7}
\end{align}
Recall the notation in \eqref{copro1} and \eqref{copro2}.
By applying the projection $\pi_{a(z,t)}$ on \eqref{complicated7}  we get the Laurent series with coefficients in $L=(L_0)^{\otimes c}$,
\begin{align}
&\lim_{z_{j,s}\to z}\bigg(\bigg(\prod_{1\leq i\leq j\leq k}\prod_{\substack{1\leq r\leq m_i\\ 1\leq s\leq m_j\\(i,r)<(j,s)}}p(z_{j,s}/z_{i,r})\bigg.\bigg)\nonumber\\
&\qquad\bigg.
\,\cdot\,\prod_{j=1}^{k} \left(z^{r_j +1} x^{(m_j)}(z_{j,1}\q^{M_j })x^{(m_j -1)}(z_{j,2}\q^{M_j +1}) ... x^{(1)}(z_{j,m_j}\q^{M_j+m_j -1})\right)\bigg)\, 1.
\label{projected}\end{align}
Formula \eqref{projected} follows from \cite[Lemma 18]{Ko1} (cf. also \cite{DF}). Roughly speaking, the projection $\pi_{a(z,t)}$
forces the operators $x(z_{j,1}\q^{M_j }),..., x(z_{j,m_j}\q^{M_j+m_j -1})$, $j=1,2,...,k$, appearing in \eqref{complicated7} to spread along the $m_j$ leftmost tensor factors of $(L_0)^{\otimes c}$ in the ordering as above, i.e.
\begin{align*}
&x(z_{j,1}\q^{M_j },t)\text{ appears on the tensor factor }m_j;\\
&x(z_{j,2}\q^{M_{j}+1 },t)\text{ appears on the tensor factor }m_j-1;\\
&\qquad\vdots\\
&x(z_{j,m_j}\q^{M_j+m_j -1},t)\text{ appears on the tensor factor }1.
\end{align*}
 For more details the reader may consult \cite[Subsection 4.1]{Ko1}.

Consider the expression  under the limit in \eqref{projected}. Formulas \eqref{copro1}, \eqref{copro2}, \eqref{r1} and \eqref{normal1} imply that its every tensor factor consists (up to a nonzero multiplicative scalar) of the operators positioned in the following order: $$\phi(\,\,) \cdots \phi(\,\,) :x(\,\,)\cdots x(\,\,):,$$
(with the appropriate arguments) and multiplied by some polynomials in
$\mathbb{C}(q^{1/2})[z_{j,s}/z_{i,r}]_0$.
Hence, the limit $\lim_{z_{j,s}\to z}$ of each tensor factor, when applied on $1\in L_0$, is   nonzero, so the statement of Lemma clearly follows.
\end{prf}

\begin{lem}\label{independent}
The set of all nonzero monomials
\begin{align}
&A_{1,l_1,m_1}^{r_1}(z,t)_{-1}A_{2,l_2,m_2}^{r_2}(z,t)_{-1}\cdots A_{k-1,l_{k-1},m_{k-1}}^{r_{k-1}}(z,t)_{-1} A_{k,l_k,m_k}^{r_k}(z,t),\label{nnzr}\\
&r_1,...,r_{k-1}\leq -2,\, r_{k}\leq -1,\, 0\leq l_j\leq -r_j -1,\, 1\leq m_j \leq c,\, j=1,...,k,\, k\in\mathbb{Z}_{\geq 0},\nonumber
\end{align}
is linearly independent.
\end{lem}

\begin{prf}
Assume that nonzero  monomials in \eqref{nnzr} are not linearly independent.
Let 
\begin{equation}\label{zero1}
\sum_{i=1}^{n} \gamma_i A_{i}(z,t) =0
\end{equation}
be a linear combination of the given nonzero monomials
\begin{align*}
A_i(z,t)=A_{1,l_{1,i},m_{1,i}}^{r_{1,i}}(z,t)_{-1}\cdots A_{k_i-1,l_{k_i-1,i},m_{k_i-1,i}}^{r_{k_i -1,i}}(z,t)_{-1} A_{k_i,l_{k_i,i},m_{k_i,i}}^{r_{k_i ,i}}(z,t),
\end{align*}
where
\begin{align*}
&r_{1,i},...,r_{k_i -1,i}\leq -2,\,\, r_{k_i,i}\leq -1,\,\, 0\leq l_{j,i}\leq -r_{j,i} -1,\,\, 1\leq m_{j,i} \leq c,\\
& j=1,2,...,k_i,\,\, k_i\in\mathbb{Z}_{\geq 0},\,\, i=1,2,...,n,
\end{align*}
such that all the scalars $\gamma_i$ are nonzero and $n>1$ is minimal.  Set
$C_i = m_{1,i}+...+ m_{k_i,i}$ for $i=1,2,...,n$. Applying the operator $\underbrace{q^{\alpha}\otimes...\otimes q^{\alpha}}_{c}$ on \eqref{zero1}  we get
\begin{equation*}
\sum_{i=1}^{n} q^{2C_i}\gamma_i A_{i}(z,t) =0.
\end{equation*}
Since $n$ was minimal, this implies $C_i =C_j$ for all $i,j=1,2,...,n$. We will denote $C_j$ by $C$. 

Set 
\begin{equation}\label{dovi}
D_i=C-r_{1,i}-...- r_{k_i,i}-k_i,\quad i=1,...,n.
\end{equation} 
Then $\deg_t A_{i}(z,t)=D_i $ for all $i=1,...,n$.
Without loss of generality we can assume that $D_1\leq D_2\leq ...\leq D_n$. Multiplying \eqref{zero1}
by $t^{-D_1}$, setting $t=0$ and then multiplying the resulting equality by $t^{D_1}$ we get another linear constraint among
$A_{i}(z,t)$. If $D_i\neq D_j$ for some $i,j=1,2,...,n$, this gives a contradiction to minimality of $n$. Therefore, we can assume
that $D_i =D_j$ for all $i,j=1,2,...,n$, i.e. $\deg_t A_{i}(z,t)=\deg_t A_{j}(z,t)$.

Every $A_i(z,t)$ can be expanded as in \eqref{zero2}:
\begin{align}
&A_i(z,t)=\lim_{z_{j,s}\to z}\bigg(\bigg(\prod_{1\leq u\leq v\leq k}\prod_{\substack{1\leq r\leq m_u\\ 1\leq s\leq m_v\\(u,r)<(v,s)}}p_i(z_{v,s}/z_{u,r})
\bigg. \bigg)\nonumber\\
&\bigg.\cdot\prod_{j=1}^{k_i} \left((z^{r_{j,i} +1}\otimes t^{-r_{j,i} -1}) x(z_{j,1}\q^{M_{j,i}+l_{j,i} },t)x(z_{j,2}\q^{M_{j,i}+l_{j,i} +1},t) ... x(z_{j,m_j}\q^{M_{j,i}+l_{j,i}+m_{j,i} -1},t)\right)\bigg),\nonumber\\
&l_{j,i}=0,...,-r_{j,i}-1, \,\, M_{j,i}=\sum_{t=1}^{j-1}(-r_{t,i} -1), \,\, j=1,..,k_i,\,\, i=1,...,n.\label{zero4}
\end{align}
We can assume that the monomials $A_i(z,t)$
are different.
 Consider the powers of $\q$ in the arguments of the operators $x(\,\,)$ in \eqref{zero4}, 
\begin{equation}\label{quovi}
\q^{M_{j,i}+l_{j,i}},\,\q^{M_{j,i}+l_{j,i} +1},...\,,\,\q^{M_{j,i}+l_{j,i}+m_{j,i} -1}.
\end{equation}
There exists an integer $p$ such that $\q^{p}$ appears in  \eqref{quovi} for $i=n$ and does not 
appear in 
\eqref{quovi} for $i=1$. Define the (invertible) operators $\mathcal{E}_c^{\pm 1} (z)$ on $(L_0)^{\otimes c}$ by
$$\mathcal{E}_c^{\pm 1} (z)=\mathcal{E}_- (z\q^p)^{\pm 1}\otimes \mathcal{E}_- (z\q^{p+1})^{\pm 1}\otimes\cdots\otimes\mathcal{E}_- (z\q^{p+c-1})^{\pm 1}.$$
Recall \eqref{komutiranjesE-} and \eqref{785}.
Applying these operators on \eqref{zero1} we get
\begin{align*}
0=\mathcal{E}_c^{-1} (z)\left(\sum_{i=1}^{n} \gamma_i A_{i}(z,t)\right)\mathcal{E}_c (z)=\sum_{i=1}^{n} \gamma_i \mathcal{E}_c^{-1} (z)A_{i}(z,t)\mathcal{E}_c (z)=\sum_{i=1}^{n} \gamma_i\beta_i A_{i}(z,t),
\end{align*}
for some scalars $\beta_i$, $i=1,2,...,n$, such that $\beta_1 \neq 0$ and $\beta_n =0$, which is in contradiction to minimality of $n$.
\end{prf}

As we observed earlier, the elements $a(z,t),b(z,t)\in \mathcal{B}_{c,1}$ can be expressed as linear combinations of monomials $\eqref{azt}$. Certain summands in such linear combinations for $a(z,t)$ and $b(z,t)$ may coincide. However, each element in $\mathcal{B}_{c,1}$ is uniquely determined by its leading term; $a(z,t)$ and $b(z,t)$ are 
equal if and only if their leading terms are equal. 
Notice that Lemmas \ref{zero5} and \ref{independent} imply that all the elements in $\mathcal{B}_{c,1}$ are nonzero.

We can introduce a linear ordering among leading terms. Let
\begin{align*}
&A(z,t)=A_{1,0,m_1}^{r_1}(z,t)_{-1}A_{1,0,m_2}^{r_2}(z,t)_{-1}\cdots A_{k-1,0,m_{k-1}}^{r_{k-1}}(z,t)_{-1} A_{k,0,m_k}^{r_k}(z,t);\\
&B(z,t)=A_{1,0,n_1}^{s_1}(z,t)_{-1}A_{1,0,n_2}^{s_2}(z,t)_{-1}\cdots A_{l-1,0,n_{l-1}}^{s_{l-1}}(z,t)_{-1} A_{l,0,n_l}^{s_l}(z,t).
\end{align*}
We write
\begin{align}
A(z,t)\prec B(z,t)
\quad\text{if}\quad&\sum_{j=1}^{k}(m_j-r_j)-k<\sum_{j=1}^{l}(n_j-s_j)-l\nonumber\\
\quad\text{or}\quad&\sum_{j=1}^{k}(m_j-r_j)-k=\sum_{j=1}^{l}(n_j-s_j)-l\nonumber\\
\quad\qquad&\text{and}\quad(m_1,r_1,m_2,r_2,...,m_k,r_k)\nonumber\\
&\qquad\qquad <(n_1,s_1,n_2,s_2,...,n_l,s_l),\label{ordering}
\end{align}
where "$<$" in \eqref{ordering} is the usual lexicographic ordering.

\begin{lem}\label{finallemma}
The set $\mathcal{B}_{c,1}$ is linearly independent.
\end{lem}

\begin{prf}
Let 
\begin{equation*}
\sum_{i=1}^{n} \gamma_i a_{i} (z,t) =0
\end{equation*}
be a linear combination of different elements $a_{i} (z,t)\in \mathcal{B}_{c,1}$ such that the scalars $\gamma_i$
are nonzero. 
Suppose that the leading term of $a_j (z,t)$ is greater, regarding the ordering "$\prec$", than the leading terms of $a_{i} (z,t)$, $i\neq j$.
When we express $a_{i} (z,t)$, $i\neq j$, as a linear combination of summands \eqref{azt}, none of these summands equals the leading term of
$a_j (z,t)$. Therefore, Lemma \ref{independent} implies $\gamma_j =0$, which is in contradiction to our choice of scalars.
\end{prf}

Finally, the main result of this section follows from Proposition \ref{spanningset7} and Lemma \ref{finallemma}:

\begin{thm}\label{tralalala}
The set $\mathcal{B}_{c,1}$ forms a basis for $V_{c,1}$.
\end{thm}

The fact that we considered  $\mathcal{E}(L)_t$ instead of $\mathcal{E}(L)$ did not affect the form of the basis $\mathcal{B}_{c,1}$, i.e. {\em the size} of the space $V_{c,1}$:

\begin{kor}\label{animportantcor}
The set $\left\{a(z,t)\left|_{t=1}\,:\,\right. a(z,t)\in\mathcal{B}_{c,1}\right\}$ is linearly independent.
\end{kor}

\begin{prf}
By checking the proofs of Lemmas \ref{zero5}, \ref{independent} and \ref{finallemma}, we see that  all the arguments except one remain valid for
$a(z,t)\left|_{t=1}\right.$.
In the proof of Lemma \ref{independent}, we considered linear combination \eqref{zero1}, $\sum_{i=1}^{n} \gamma_i A_{i}(z,t) =0$ and showed that we can assume that all the elements $D_i$,
$i=1,2,...,n$, defined by \eqref{dovi}, are equal.
Even though such a conclusion does not hold for $t=1$, we can always derive from \eqref{zero1} a new linear combination,
$\sum_{j=1}^{m} \beta_j (z) B_{j}(z,t) =0$, where $\beta_j(z)$ are nonzero polynomials in $z^{-1}$, $B_{j}(z,t)\neq B_{k}(z,t)z^l$ for all $j,k=1,2,...,m$, $l\in\mathbb{Z}$, $j\neq k$ and $B_{j}(z,t)\in\left\{A_{1}(z,t),A_{2}(z,t),...,A_{n}(z,t)\right\}$ for all $j=1,...,m$. Note that $ \beta_j (z) B_{j}(z,t)\neq 0$ for all $j=1,...,m$. Now, with this linear combination, we can proceed in the same way as in the proof of Lemma \ref{independent}.
\end{prf}

Corollary \ref{animportantcor} implies that the restriction of the evaluation homomorphism $$\varphi\colon a(z)\otimes g(t)\mapsto a(z)\otimes g(1)$$
on $V_{c,1}$ is injective, so 
$\varphi\colon V_{c,1}\to V_c$ is an isomorphism of vector spaces $V_{c,1}$ and $V_c=\im \varphi$.  Hence, we can define the structure
of graded nonlocal $\q$-vertex algebra on $V_c \subset \mathcal{E}(L)$ in the following way:
\begin{align}
&Y(a(z),z_0)b(z)=\sum_{r\in\mathbb{Z}} \varphi(\varphi^{-1}(a(z))_r \varphi^{-1}(b(z))) z_{0}^{-r-1};\label{fgh11}\\
&\vac =1_{\mathcal{E}(L)};\label{fgh22}\\
&(V_c)_{(s)} = \varphi((V_{c,1})_{(s)})\label{fgh33}
\end{align}
for all $a(z),b(z)\in V_c$, $s\in\mathbb{Z}$. Note that $x(z)$ is an element of  $V_c$, so, in view of \eqref{fgh11}--\eqref{fgh33},  we can say that $V_c$ is the graded nonlocal $\q$-vertex algebra
generated by $x(z)$.

\begin{rem}
It is not clear whether the techniques, which we used to prove the linear independence of the set $\mathcal{B}_{c,n}$ for $c\geq 1$, $n=1$, can be applied to the case
 $c>1$, 
$n=0$. The similar spanning set for $V_{c,0}$
 can be easily constructed, but  some of its elements have zero leading terms 
and, furthermore, some different elements  have equal (nonzero) leading terms. 

Recall the notation from Definition \ref{definition}. The case $c=1$, $n=0$ was solved in \cite{Ko4} using the slightly different vertex algebraic setting. The results therein imply that the set 
\begin{align*}
\mathcal{B}_{1,0}=\big\{&x_{1,0}(z,t)_{r_1}x_{1,0}(z,t)_{r_2}\ldots x_{1,0}(z,t)_{r_k}\vac\,:\big.\big. r_1,...,r_{k-1}\leq -3,\,r_k\leq -1,\,k\in\mathbb{Z}_{\geq 0} \big\}
\end{align*}
forms a basis for $V_{1,0}$ and that the analogue of Corollary \ref{animportantcor} holds.
Moreover, the leading terms of $\varphi(\mathcal{B}_{1,0})$ and $\varphi(\mathcal{B}_{1,1})$
coincide up to a nonzero multiplicative factor.

On the other hand, construction of similar bases
for $c\geq 1$, $n>1$ is straightforward since, due to the definition of the vertex operator products, quantum integrability relations \eqref{ding_miwa} are not applicable anymore.
The set
\begin{align*}
\mathcal{B}_{c,n}=\big\{&x_{m_1,n}(z,t)_{r_1}x_{m_2,n}(z,t)_{r_2}\ldots x_{m_k,n}(z,t)_{r_k}\vac\,:\big.\\
&\big. r_1,...,r_{k-1}\leq -2,\,r_k\leq -1,\, m_j\geq 1,\,j=1,...,k,\,k\in\mathbb{Z}_{\geq 0} \big\}
\end{align*}
forms a basis of $V_{c,n}$ and its
linear independence  can be proved by
using the arguments from Lemma \ref{independent}.
\end{rem}

\subsection{Character formulae}

Let
\begin{equation}\label{diff2}
a(z,t)=x_{m_1}(z,t)_{r_1}x_{m_2}(z,t)_{r_2}\ldots x_{m_k}(z,t)_{r_k}\vac\in \mathcal{B}_{c,1}
\end{equation}
be an arbitrary basis element.
We can associate a colored Young diagram with $a(z,t)$ in the following way.  The diagram will have $k$ columns: the first column will have 
$-(r_1 + ... + r_k)$ rows, the second column will have $-(r_2 + ... + r_k)$ rows, ..., the last column will have $-r_k$ rows. 
Boxes in column $j=1,2,...,k$ will be labeled by $m_j$.

\begin{figure}
\begin{tikzpicture}[scale=0.8]
\tikzstyle{every node}=[font=\tiny]
\node at (8.5,0.5) {$7$};
\node at (8.5,1.5) {$7$};
\node at (8.5,2.5) {$7$};
\node at (8.5,3.5) {$7$};
\node at (8.5,4.5) {$7$};
\node at (8.5,5.5) {$7$};
\node at (8.5,6.5) {$7$};
\node at (8.5,7.5) {$7$};
\node at (9.5,3.5) {$4$};
\node at (9.5,4.5) {$4$};
\node at (9.5,5.5) {$4$};
\node at (9.5,6.5) {$4$};
\node at (9.5,7.5) {$4$};
\node at (10.5,5.5) {$5$};
\node at (10.5,6.5) {$5$};
\node at (10.5,7.5) {$5$};
\node at (11.5,7.5) {$6$};
\draw  (8,1) -- (8,0) -- (9,0) -- (9,1) -- (8,1);
\draw  (8,2) -- (8,1) -- (9,1) -- (9,2) -- (8,2);
\draw  (8,3) -- (8,2) -- (9,2) -- (9,3) -- (8,3);
\draw  (8,4) -- (8,3) -- (9,3) -- (9,4) -- (8,4);
\draw  (8,5) -- (8,4) -- (9,4) -- (9,5) -- (8,5);
\draw  (8,6) -- (8,5) -- (9,5) -- (9,6) -- (8,6);
\draw  (8,7) -- (8,6) -- (9,6) -- (9,7) -- (8,7);
\draw  (8,8) -- (8,7) -- (9,7) -- (9,8) -- (8,8);

\draw  (9,4) -- (9,3) -- (10,3) -- (10,4) -- (9,4);
\draw  (9,5) -- (9,4) -- (10,4) -- (10,5) -- (9,5);
\draw  (9,6) -- (9,5) -- (10,5) -- (10,6) -- (9,6);
\draw  (9,7) -- (9,6) -- (10,6) -- (10,7) -- (9,7);
\draw  (9,8) -- (9,7) -- (10,7) -- (10,8) -- (9,8);

\draw  (10,6) -- (10,5) -- (11,5) -- (11,6) -- (10,6);
\draw  (10,7) -- (10,6) -- (11,6) -- (11,7) -- (10,7);
\draw  (10,8) -- (10,7) -- (11,7) -- (11,8) -- (10,8);

\draw  (11,8) -- (11,7) -- (12,7) -- (12,8) -- (11,8);

\end{tikzpicture}\caption{Diagram for $b(z,t)=x_{7}(z,t)_{-3}x_{4}(z,t)_{-2} x_{5}(z,t)_{-2}x_{6}(z,t)_{-1}\vac$}
\label{pic4}
\end{figure}

For the  basis element $a(z,t)$
in \eqref{diff2}
set
\begin{align}\label{qdeg}
&\deg_\q a(z,t) =-(r_1 +2r_2 +3r_3+...+kr_k).
\end{align}
Note that $\deg_\q a(z,t)$ is equal to the number of boxes in the diagram for $a(z,t)$. Next, define $\deg_\q \vac =0$.

\begin{ex}
The diagram for the basis element
$$b(z,t)=x_{7}(z,t)_{-3}x_{4}(z,t)_{-2} x_{5}(z,t)_{-2}x_{6}(z,t)_{-1}\vac\in \mathcal{B}_{c,1}$$
is given in Figure \ref{pic4}. Of course,  we assume here that $c\geq 7$. The degree $\deg_\q b(z,t)$ equals $17$.
\end{ex}

The space $V_{c,1}$ admits a decomposition
$$V_{c,1}=\coprod_{r\geq 0} (V_{c,1})_{r},\quad (V_{c,1})_{r}=\left\{a(z,t)\in V_{c,1}\,:\,\deg_\q a(z,t)=r\right\},$$
and all the subspaces $(V_{c,1})_{r}$ are finite-dimensional. Hence, we can define
$$\textstyle\ch_{\q} V_{c,1} =\displaystyle\sum_{r\geq 0} \dim (V_{c,1})_{r}\, q^r.$$
Naturally, parameters $\q$ and $q$ in the above formula are not related.

For the basis element $a(z,t)$ in \eqref{diff2} and $j=1,2,...,k$ set
\begin{align*}
&a_j (z,t)= x_{m_j}(z,t)_{r_j}x_{m_{j+1}}(z,t)_{r_{j+1}}\ldots x_{m_k}(z,t)_{r_k}\vac\in \mathcal{B}_{c,1};\\
&a_{k+1}(z,t)=\vac\in \mathcal{B}_{c,1}.
\end{align*}
Notice that the defining conditions for the basis $\mathcal{B}_{c,1}$ in \eqref{basis785}, 
\begin{align}
&r_j\leq -2,\quad\text{for all } j=1,2,...,k-1;\label{d2c1}\\
& r_{k}\leq -1\label{d2c2}
\end{align}
can be expressed in terms of $\deg_\q$  by
\begin{align}
&\left(\deg_\q a_{j+1} (z,t) -\deg_\q a_{j+2} (z,t)\right)-\left(\deg_\q a_{j} (z,t) -\deg_\q a_{j+1} (z,t)\right) \leq -2\label{d2c3}\\
&\qquad \text{for all }j=1,...,k-1;\nonumber\\
&\deg_\q a_{k}(z,t)\leq -1.\label{d2c4}
\end{align}

\begin{thm}\label{charformula}
\begin{equation}\label{qrr}
\textstyle\ch_{\q} V_{c,1} =\displaystyle\sum_{r\geq 0} \frac{q^{r^2}}{(1-q)(1-q^2)\cdots(1-q^r)}c^r.
\end{equation}
\end{thm}

\begin{prf}
Let $c=1$.
Coefficient of  $q^n$ in the $r$th summand $q^{r^2}(1-q)^{-1}\cdots(1-q^r)^{-1}$
equals  the number of partitions of $n$ into $r$ parts such that the difference among consecutive parts is at least two.
Therefore, equivalence of \eqref{d2c1}, \eqref{d2c2} and \eqref{d2c3}, \eqref{d2c4} implies that the character formula \eqref{qrr} holds for $c=1$. For an arbitrary level $c$, the basis $\mathcal{B}_{c,1}$
consists of monomials of quasi-particles, which have charges $1,2,...,c$ and satisfy the same  conditions \eqref{d2c3} and  \eqref{d2c4}, while the  
factor $c^r$ in \eqref{qrr} counts all the possible choices of charges in a monomial consisting of $r$
quasi-particles.
\end{prf}

\begin{rem}
In contrast to the function $\deg_\q$, defined  in \eqref{qdeg}, the original grading  for $V_{c,1}$, which is given by the function  $\wt$ in \eqref{newlabel},
depends on the labels $m_1,...,m_k$ of the element $a(z,t)$ in \eqref{diff2}.
Defining conditions \eqref{d2c1} and \eqref{d2c2} can be expressed in terms of $\wt$ by
\begin{align*}
& \wt a_{j} (z,t) - \wt a_{j+1} (z,t)\geq m_{j}+1
\quad \text{for all }j=1,...,k-1;\\
&\wt a_{k}(z,t)\geq 1.
\end{align*}
Furthermore, the value $\wt a(z,t)$ can be easily obtained from the corresponding diagram:
$$\wt a(z,t)\,=\,\text{(number of rows)}\,-\,\text{(number of columns)}\,+\,\text{(sum of labels)}.$$
For example (see Figure \ref{pic4}),
$$\wt b(z,t)\, =\, 8\, -\, 4\,+\,(7+4+5+6)\,=\, 26.$$
\end{rem}

\section*{Acknowledgement}
The research was supported by the Australian Research Council
and  by the Croatian Science Foundation under the project 2634.
The author would like to thank Mirko Primc for comments that greatly improved the manuscript.

\end{document}